\documentclass[a4paper,reqno,12pt]{amsart}
\usepackage{amsmath,amssymb,amsthm,mathrsfs,mathtools,amsfonts}
\usepackage[left=25mm,right=25mm,top=24mm,bottom=25mm,headsep=6mm]{geometry}
\usepackage[colorlinks=true,linkcolor=blue,citecolor=blue,urlcolor=cyan]{hyperref}
\hypersetup{pdftitle={Parabolic star products and the relative braid restriction}}
\numberwithin{equation}{section}
\allowdisplaybreaks[1]
\AtBeginDocument{\setlength{\abovedisplayskip}{7pt plus 2pt minus 2pt}
\setlength{\belowdisplayskip}{7pt plus 2pt minus 2pt}}

\newtheorem{proposition}{Proposition}[section]
\newtheorem{lemma}[proposition]{Lemma}
\newtheorem{theorem}[proposition]{Theorem}
\newtheorem{corollary}[proposition]{Corollary}
\newtheorem{definition}[proposition]{Definition}
\theoremstyle{definition}
\newtheorem{remark}[proposition]{Remark}

\theoremstyle{plain}
\newtheorem{introtheorem}{Theorem}

\theoremstyle{definition}

\def\fg{\mathfrak g}
\def\Q{\mathbb Q}
\def\K{\mathbb K}

\def\N{\mathbb N}
\def\Z{\mathbb Z}
\def\I{\mathbb I}
\def\J{\mathbb J}
\def \ba{{\mathbf a}}
\def \bi{{\mathbf i}}
\def \bj{{\mathbf j}}
\def \bk{{\mathbf k}}
\def\bI{\I_\bullet}
\def\wI{\I_\circ}
\def\bIi{\I_{\bullet,i}}
\def\bw{w_\bullet}
\def\bwi{w_{\bullet,i}}
\def\br{\mathbf r}
\def\tU{\widetilde{\mathbf U}}

\def\tUi{\widetilde{\mathbf U}^{\imath}}

\def \tbU{\tU_{\bullet}}
\def \tbiU{\tU_{\bullet,i}}

\def\tTD{\widetilde T}
\def\tT{\widetilde{\mathscr T}}
\def\tTT{\widetilde{\mathbf T}}
\def\tPsi{\widetilde\Psi}
\def\tk{\widetilde k}
\def\tfX{\widetilde\Upsilon}

\def\cA{\mathcal A}
\def\cF{\mathcal F}
\def\cR{\mathcal R}
\def\cH{\mathcal H}
\def\cM{\mathcal M}

\def\cJ{\mathcal J}
\def\cV{\mathcal V}
\def\cK{\mathcal K}
\def\cD{\mathcal D}
\def\cG{\mathcal G}
\def\vs{\varsigma}
\def\bvs{\boldsymbol\varsigma}
\def\dm{\diamond}
\def \ov{\overline}
\def \bvs{{\boldsymbol{\varsigma}}}
\def \wItau{\I_{\circ,\tau}}
\def \fwItau{\wItau^{\mathrm{fin}}}
\def \Ob{\Omega_\bullet}
\def \wObb{\widetilde{\mathbf{\Omega}}_\bullet}
\def \Gb{\Gamma_\bullet}

\newcommand{\tTa}[1]{\tTT'_{#1,-1}}
\newcommand{\id}{\operatorname{id}}
\newcommand{\ad}{\operatorname{ad}}
\newcommand{\Span}{\operatorname{span}}
\newcommand{\eps}{\varepsilon}
\newcommand{\poly}{\mathrm{poly}}
\newcommand{\bDelta}{\overline\Delta}
\newcommand{\diag}{\mathrm{diag}}
\newcommand{\qbinom}[2]{\begin{bmatrix} #1\\#2 \end{bmatrix} }
\newcommand{\Aut}{\operatorname{Aut}\nolimits}
\newcommand{\obar}{\overline{\phantom{o}}}
\newcommand{\gr}{\operatorname{gr}}
\newcommand{\wt}{\operatorname{wt}}
\newcommand{\ev}{\operatorname{ev}}

	\newcommand{\redtext}[1]{\textcolor{red}{#1}}

	\newcommand{\red}[1]{\redtext{ #1}}

\begin{document}	
%%%%%%%%%
\title[$\mathrm i$Braid group symmetries of $\mathrm i$quantum groups]{$\mathrm i$Braid group symmetries of $\mathrm i$quantum groups via star product}

	\author[Ming Lu]{Ming Lu}
	\address{Department of Mathematics, Sichuan University, Chengdu 610064, P.R. China}
	\email{luming@scu.edu.cn}

\author[Xujie Zhou]{Xujie Zhou}
	\address{Department of Mathematics, Sichuan University, Chengdu 610064, P.R. China}
\email {2025222010028@stu.scu.edu.cn}

	
\begin{abstract}
We establish ibraid group symmetries on universal iquantum groups of symmetrizable Kac–Moody type, proving a conjecture of Wang and Zhang. Building on the star product realization of quantum symmetric pairs $(\tU,\tUi)$, we relate these symmetries to Lusztig's braid group symmetries, extending the compatibility previously obtained via iHopf algebras in the quasi-split setting. Via the Letzter map $\psi$ from the iquantum group $\tUi$ to the horospherical algebra $\cA\subset \tU$, we show that the ibraid group symmetry $\tTa{i}$ of $\tUi$ transported to  \(\mathcal A\) coincides with Lusztig's braid group symmetry \(\widetilde T_{\mathbf r_i}^{-1}\) on the subalgebra \(\mathcal R_{\bullet,i}^{-,r}\). Finally, a mixed product formula, together with a parabolic decomposition, yields an explicit formula for the transported action on the entire horospherical algebra $\mathcal A$.
\end{abstract}

%We study relative braid group symmetries of universal $\mathrm i$quantum groups through the star product realization of quantum symmetric pairs.  Let $\psi:\tUi\rightarrow\cA$ be the Letzter map, which identifies $\tUi$ linearly with the horospherical algebra $\cA$ and transports the multiplication of $\tUi$ to the star product on $\cA$.  For every white node $i$, we show that the relative braid symmetry transported to $\cA$ restricts on the parabolic factor $\cR_{\bullet,i}^{-,r}$ to the ordinary Lusztig's braid group symmetry $\tTD_{\br_i}^{-1}$.  We further derive an explicit formula for the transported relative braid group symmetries on the full horospherical algebra $\cA$.  

	\maketitle
    \setcounter{tocdepth}{1}
	\tableofcontents

%%%%%%%
\section{Introduction}

Quantum symmetric pairs $(\mathbf{U},\mathbf{U}^\imath)$, consisting of Drinfeld--Jimbo quantum groups $\mathbf{U}$ and their coideal subalgebras $\mathbf{U}^\imath$, were introduced by Letzter \cite{Let99,Let02} and developed systematically by Kolb \cite{Ko14}. These subalgebras $\mathbf{U}^\imath$, attached to  admissible Satake diagrams $(\I=\wI\sqcup\bI,\tau)$, are called $\mathrm i$quantum groups (also called quantum symmetric pair coideal subalgebras). They are a vast generalization of quantum groups, since quantum groups can be viewed as iquantum groups of diagonal type. %The iprogram, initiated by Bao and Wang \cite{BW18a}, is to generalize fundamental constructions of quantum groups to iquantum groups.  
The theory of their bar involutions, canonical bases, braid group symmetries, and universal $K$-matrices has since been developed in works including \cite{KP11,BK15,BW18a,BW18b,BK19,AV22,LW21b,WZ23,Z23}. 

Inspired by the iHall algebra realization \cite{LW22a}, the universal iquantum groups $\tUi$ were introduced as subalgebras of Drinfeld double quantum groups $\tU$; see \cite{WZ23}. The classical iquantum groups $\mathbf{U}^\imath$ can be recovered from $\tUi$ by central reduction. We focus on $\tUi$ in this paper.  
 %In the universal form used in this paper, an $\mathrm i$quantum group $\tUi$ is attached to an admissible Satake diagram $(\I=\wI\sqcup\bI,\tau)$ and sits inside the Drinfeld double quantum group $\tU$ as a right coideal subalgebra.

Braid group symmetries are among the fundamental structural features of quantum groups.  Lusztig constructed braid group automorphisms associated with the simple reflections of the Weyl group \cite{Lus90,Lus93}.  For quantum symmetric pairs, the corresponding symmetries are governed by the relative Weyl group.  Early formulas and conjectures appeared in \cite{KP11}; the quasi-split case was realized by reflection functors of iquivers in \cite{LW21b}, and the intrinsic construction for finite type and locally quasi-split Kac-Moody type was established in \cite{WZ23,Z23} via $K$-matrices.  We refer to these automorphisms as relative braid group (or ibraid group) symmetries.  Although the ordinary and relative braid group actions are defined on different algebras, the close relation between them is a recurring feature of the subject.

This topic was studied in \cite{CLPRW25} via iHopf algebras.  In \emph{loc. cit.}, an $\mathrm i$Hopf algebra structure was introduced on a Hopf algebra equipped with a Hopf pairing, and the resulting twisted multiplication on a Borel quantum group was shown to realize quasi-split universal $\mathrm i$quantum groups.  One application is a direct connection between Lusztig's braid group symmetries and the ibraid  group symmetries in the quasi-split setting.  On suitable Borel root vectors, the ibraid group symmetry is identified with the ordinary braid group symmetry associated with the corresponding relative simple reflection.
This was strengthened to an explicit formula for the transported ibraid group symmetry on the full Borel subalgebra; see \cite[Theorem 3.11]{CLPRW26}.

On the other hand, 
Kolb and Yakimov introduced a star product construction for quantum symmetric pairs associated with Nichols algebras of diagonal type \cite{KY20} and used it to study defining relations of $\mathbf{U}^\imath$ in \cite{KY21}.  Their recent work \cite{KY26} develops short star products for general quantum symmetric pairs. The $\mathrm i$Hopf construction and the Kolb--Yakimov's star product construction are quite different in origin, although they lead to closely related deformation formulas.  The present paper provides a star product counterpart of these compatibility results, allowing an arbitrary black part $\bI$.

The main results of this paper can be summarized in the following
theorems.

Let $(\I=\wI\cup\bI,\tau)$ be an admissible Satake diagram of
Kac--Moody type, and let $\tUi\subseteq\tU$ be the corresponding
universal iquantum group. Fix a set $\wItau$ of representatives of
the $\tau$-orbits in $\wI$, and set
\[
  \bIi=\bI\cup\{i,\tau i\},\qquad
  \fwItau=\{i\in\wItau\mid\bIi\text{ is of finite type}\}.
\]
Let $\bw$ be the longest element of $W_\bullet$ and put
$\theta=-\bw\circ\tau$. For $i\in\fwItau$, let $\bwi$ be the
longest element of $W_{\bullet,i}$ and set
\[
  \br_i=\bwi\bw,\qquad
  W^\circ=\langle\br_i\mid i\in\fwItau\rangle.
\]
Write $\tfX_i$ for the rank-one quasi $K$-matrix associated with
$\bIi$ and $\tT_{\br_i}^{\pm1}$ for the rescaled Lusztig's braid group
symmetries defined in \eqref{eq: T_i twisted def}. 

\begin{introtheorem}[Theorem~\ref{thm: relative T_i conjugate}]
\label{thm:intro-existence}
For every $i\in\fwItau$, there are mutually inverse algebra
automorphisms $\tTT'_{i,-1}$ and $\tTT''_{i,+1}$ of $\tUi$ such that,
for $x\in\tUi$,
\begin{align*}
  \tTT'_{i,-1}(x)\tfX_i
    &=\tfX_i\tT_{\br_i}^{-1}(x),\\
  \tTT''_{i,+1}(x)\tT_{\br_i}(\tfX_i)^{-1}
    &=\tT_{\br_i}(\tfX_i)^{-1}\tT_{\br_i}(x).
\end{align*}
The automorphisms $\tTT''_{i,+1}$, $i\in\fwItau$, satisfy the braid
relations of $W^\circ$.
\end{introtheorem}
The key ingredient in the proof of Theorem~\ref{thm:intro-existence}
is the bigraded subalgebra $\tU_{\text{int}}$ in
\eqref{eq:sum of bigraded subspace and int subalge}, containing
$\tUi$ and $\{F_j,\tTD_{\bw}(E_{\tau j})K_j'\mid j\in\wI\}$, together with the algebra endomorphism
$\Xi:\tU_{\text{int}}\longrightarrow\tU_{\text{int}}$ defined in
\eqref{eq:formula of Xi}. Write $\cK_\lambda=K_\lambda K_{-\theta\lambda}'$. 
This subalgebra is stable under $\tT_{\br_i}^{\pm1}$ and the anti-involution $\sigma$. The endomorphism $\Xi$ also fixes $\tbU$ and $\cK_j$ from \eqref{eq:vkHla}.
For $i\in\fwItau$ and $j\in\wI$, by Lemmas~\ref{lem:Xi is alge-homomorphism}--\ref{lem:Xi commute brT}, we obtain
\begin{gather*}
  \Xi(B_j^\sigma)=B_j,\qquad
  \Xi(F_j)=\tTD_{\bw}(E_{\tau j})K_j',\\
  \Xi\circ\tT_{\br_i}^{\pm1}
    =\tT_{\br_i}^{\pm1}\circ\Xi.
\end{gather*}
Here $B_j^\sigma=\sigma(B_j)$.
For $j\notin\bIi$, applying $\Xi$ to the polynomial expressions
for $\tT_{\br_i}^{\pm1}(F_j)$ obtained from the right adjoint action
yields the identities conjectured in
\cite[Conjecture~5.13]{WZ23}.
The quasi $K$-matrix intertwining relations then yield
Theorem~\ref{thm:intro-existence}.

The horospherical subalgebra and the Letzter map are given by
\[
  \cA=\K\langle F_j,E_k,\tU^{\imath 0}
           \mid j\in\I,\ k\in\bI\rangle,
  \qquad
  \psi:\tUi\longrightarrow\cA.
\]
The map $\psi$ is a filtered linear isomorphism and transports the
multiplication of $\tUi$ to the star product
\[
  a\ast b=\psi\bigl(\psi^{-1}(a)\psi^{-1}(b)\bigr),
  \qquad a,b\in\cA.
\]
We first establish a tensor form of the relative braid intertwining
relation. Let $\Theta_{\bullet,i}^{\imath}$ be the tensor quasi
$K$-matrix associated with $\bIi$, and let $\overline\Theta_\bullet$
be the inverse quasi $R$-matrix for the quantum subalgebra $\tbU$
associated with $\bI$.

\begin{introtheorem}[Theorem~\ref{thm:tensor}]
\label{thm:intro-tensor}
For any $x\in\tUi$, the following identity holds in
$\tU\widehat\otimes\tU$:
\[
  \Delta\bigl(\tTa{i}(x)\bigr)
  \Theta_{\bullet,i}^{\imath}\overline\Theta_\bullet
  =\Theta_{\bullet,i}^{\imath}\overline\Theta_\bullet
   (\tTa{i}\otimes\tTD'_{\br_i,-1})\Delta(x).
\]
\end{introtheorem}

Its proof combines a factorization of the tensor quasi $K$-matrix with
the compatibility of Lusztig's braid group symmetries $\tTD_{i,-1}'$ with the coproduct.
This identity relates the two braid group actions through the
coideal structure of $\tUi$.

We next compare the two braid group actions on the right coideal subalgebra
\[
  \cR_{\bullet,i}^{-,r}
  =\tU^-\cap\tTD_{\bwi}(\tU^-).
\]
%We then have the following.

\begin{introtheorem}[Theorem~\ref{thm:restriction}]
\label{thm:intro-restriction}
For any $z\in\cR_{\bullet,i}^{-,r}$, we have
\[
  (\psi\circ\tTa{i}\circ\psi^{-1})(z)
  =\tTD_{\br_i}^{-1}(z)\in\tU^-.
\]
\end{introtheorem}

The proof uses the relation between the Letzter map and the tensor quasi
$K$-matrix to control the Cartan support of $\psi^{-1}(z)$. We introduce
a grading adapted to $\bIi$ and use a specialization at $t=0$ to
recover $z$ from $\psi^{-1}(z)$. Combining this specialization with the
quasi $K$-matrix intertwining relation gives
Theorem~\ref{thm:intro-restriction}. The same grading also determines
the Cartan support of the lower-degree terms of the star product; see
Proposition~\ref{prop:star}.

\begin{samepage}
As an application of the theorem, we recover the compatibility
of ibraid group symmetries with the relative simple reflections from
\cite[Theorem~7.13]{WZ23}. Namely, if $w\in W^\circ$ and $i\in\wI$
satisfy $wi\in\wI$, then
\[
  \tTT'_{\underline{w},-1}(B_i)=B_{wi}
\]
for some reduced expression $\underline{w}$ of $w$; see
Corollary~\ref{thm:fact1}.
\end{samepage}

To describe the action on the whole of $\cA$, set
\[
  \cH_\bullet=\tbU\tU^{\imath 0},
  \qquad
  \cR_i=\tbiU^-\cap\tTD_{\bw}(\tbiU^-).
\]
Multiplication induces the parabolic decomposition
\[
  \cA\cong
  \cH_\bullet\otimes\tTD_{\bw}^{-1}(\cR_i)
  \otimes\cR_{\bullet,i}^{-,r}.
\]
For $u\in\bigl(\tTD_{\bw}^{-1}(\cR_i)\bigr)_{-\mu}$ and
$x\in(\cR_{\bullet,i}^{-,r})_{-\nu}$, we establish the mixed product
formula in Theorem~\ref{thm:mixed}:
\[
  (\psi\circ\tTa{i}\circ\psi^{-1})(ux)
  =q^{(\mu,\nu)}\tTD_{\br_i}^{-1}(x)
    (\psi\circ\tTa{i}\circ\psi^{-1})(u).
\]
We obtain this formula using the opposite Letzter projection and the
coproduct of $u$. The remaining local action on $\tbiU^-$ is computed
explicitly in Proposition~\ref{prop:local formula for U-}. Together with
the known action on $\cH_\bullet$, these results give the following
formula on the whole of $\cA$.

\begin{samepage}
Let $\tau_{\bullet,i}$ be the
diagram involution of $\bIi$ determined by
$\bwi(\alpha_j)=-\alpha_{\tau_{\bullet,i}j}$ for $j\in\bIi$.
We use the distinguished parameter $\bvs_\dm$ and its monomials
$\vs_{\mu,\dm}$ from \eqref{eq:disting-para} and \eqref{eq:vkHla}.

\begin{introtheorem}[Theorem~\ref{thm:full-action}]
\label{thm:intro-full-action}
For $h\in\cH_\bullet$,
$u\in\bigl(\tTD_{\bw}^{-1}(\cR_i)\bigr)_{-\mu}$, and
$x\in(\cR_{\bullet,i}^{-,r})_{-\nu}$, we have
\begin{align*}
  (\psi\circ\tTa{i}\circ\psi^{-1})(hux)
  &={\vs_{\mu,\dm}}
    q^{(\mu,\mu+\theta\mu+2\nu)/2}
    \tT_{\br_i}^{-1}(h)\cdot\tTD_{\br_i}^{-1}(x)
\cdot
    \tT_{\bw}^{2}\bigl(\tau_{\bullet,i}(\tau(u))\bigr)
    \cK_{\tau_{\bullet,i}\tau\mu}^{-1}.
\end{align*}
\end{introtheorem}
\end{samepage}

All products in Theorem~\ref{thm:intro-full-action} are ordinary products in $\tU$. The Cartan datum is allowed to be of arbitrary symmetrizable Kac--Moody type; the finite-type assumptions concern $\bI$ and the subdiagrams $\bIi$ defining the relative simple reflections.

The paper is organized as follows.  Section~\ref{sec:QG-iQG}  reviews quantum groups, universal $\mathrm i$quantum groups, the Letzter map, and the star product. Section~\ref{sec:i braid action} establishes the relative braid (or ibraid) group symmetries on universal $\mathrm i$quantum groups. Section~\ref{sec:tensor K} develops the tensor quasi $K$-matrix and proves the braid--coproduct intertwining identity.  Section~\ref{sec:R} studies the parabolic factor $\cR_{\bullet,i}^{-,r}$, establishes the Cartan-support properties of the star product, and proves the restriction theorem.  Section~\ref{sec:global} introduces the opposite Letzter projection, proves the mixed product formula, and derives the explicit formula for the transported ibraid group symmetries on $\cA$.

\medskip
\noindent\textbf{Acknowledgment.}
ML is partially supported by the National Natural Science Foundation of
China (No.~12671049, 12631002).

\medskip
\noindent\textbf{Declaration on the Use of Generative AI.}
The authors used OpenAI's ChatGPT to improve the language and presentation of
the manuscript and as an additional tool to verify the logical consistency
of some proofs. All mathematical arguments, statements, and references were
independently verified by the authors, who take full responsibility for the
content of the paper.

%%%%%%%%%%%%%%%%%%%
\section{Quantum groups, iquantum groups and star product}
\label{sec:QG-iQG}
This preliminary section recalls quantum groups and iquantum groups, together with their realization via the Letzter map and the star product as in \cite{KY21,KY26}.
%%%%%%%%%%%%%
\subsection{Quantum groups}
Let $\I=\{1,\dots,n\}$. 
	Let $A=(a_{ij})_{i,j \in \I}$ be the symmetrizable generalized Cartan matrix (GCM) of a Kac--Moody Lie algebra $\fg$. Let $D=\diag(\epsilon_i\mid i\in \I)$ with $\epsilon_i\in\Z_{>0}$ be the symmetrizer of $A$, i.e., $DA$ is symmetric. 
	Let $\{\alpha_i\mid i\in\I\}$ be a set of simple roots of $\fg$, and denote the root lattice by $\Z{\I}:=\Z\alpha_1\oplus\cdots\oplus\Z\alpha_n$. Equip $\Z\I$ with the symmetric bilinear form determined by
	\begin{align} \label{BilForm}
		(\alpha_i,\alpha_j)=\epsilon_ia_{ij},\quad \forall i,j\in\I.
	\end{align}
	
	The simple reflection $s_i:\Z{\I}\rightarrow\Z{\I}$ is defined by $s_i(\alpha_j)=\alpha_j-a_{ij}\alpha_i$, for $i,j\in \I$.
	Denote the Weyl group by $W =\langle s_i\mid i\in \I\rangle$.
	
	Let $q$ be an indeterminate and $\Q(q)$  be the field of rational functions in $q$ with coefficients in $\Q$. Let $\K$ be the algebraic closure of $\Q(q)$. Set 
	$$q_i=q^{\epsilon_i},\qquad \forall i\in\I.$$
    For $r\in\N$ and $m \in \Z$, set
	\begin{align*}
	  [r]_{i}=\frac{q_i^r-q_i^{-r}}{q_i-q_i^{-1}},
	\quad
	[r]_{i}^!=\prod_{n=1}^r [n]_{i}, \quad \qbinom{m}{r}_{i} =\frac{[m]_{i}[m-1]_{i}\ldots [m-r+1]_{i}}{[r]_{i}^!}.  
	\end{align*}
    For $A,B$ in a $\K$-algebra, write $[A, B]=AB-BA$ and $A^{(n)}=\frac{A^n}{[n]_i^!}$ for the divided power.
	
	The quantum group $\tU := \tU_q(\fg)$ is defined to be the $\K$-algebra generated by $E_i,F_i, K_i,K_i'$, $i\in \I$, where $K_i,K_i'$ are invertible, subject to the following relations:  for $i, j \in \I$,
	\begin{align}
		[E_i,F_j]= \delta_{ij}\frac{K_i-K_i'}{q_i-q_i^{-1}},  &\qquad [K_i,K_j]=[K_i,K_j']  =[K_i',K_j']=0,
		\label{eq:KK}
		\\
		K_i E_j=q_i^{a_{ij}} E_j K_i, & \qquad K_i F_j=q_i^{-a_{ij}} F_j K_i,
		\label{eq:EK}
		\\
		K_i' E_j=q_i^{-a_{ij}} E_j K_i', & \qquad K_i' F_j=q_i^{a_{ij}} F_j K_i',
		\label{eq:K2}
	\end{align}
	and for $i\neq j \in \I$,
	\begin{align}
		& \sum_{r=0}^{1-a_{ij}} (-1)^r   E_i^{(r)} E_j  E_i^{(1-a_{ij}-r)}=0,\qquad \sum_{r=0}^{1-a_{ij}} (-1)^r   F_i^{(r)} F_j  F_i^{(1-a_{ij}-r)}=0.
		\label{eq:serre}
	\end{align}
	
	The algebra $\tU$ is $\Z\I$-graded by 
	\begin{equation}\label{eq:wt-grade}
	   \deg E_i=\alpha_i, \qquad\deg F_i=-\alpha_i,\qquad \deg K_i=0=\deg K_i'. 
	\end{equation} 
	Let $\tU_\mu$ denote the homogeneous subspace of degree $\mu$, so that $\tU=\oplus_{\mu\in\Z\I} \tU_\mu$. We usually call \eqref{eq:wt-grade} the $\wt$-grading and write $\wt(x)=\mu$ for $x\in\tU_\mu$.

    By \eqref{eq:KK}-\eqref{eq:serre}, the algebra $\tU$ is also $\Z\I$-graded by
    \begin{align}
      \deg_\I(E_j)=\alpha_j,\quad\deg_\I(F_j)=0,\quad\deg_\I K_j^{\pm1}=\deg_\I K_j'^{\pm1}=\pm\alpha_j,\qquad j\in\I,\label{eq:ZI-grade-1}\\
      \deg'_\I(E_j)=0,\quad\deg'_\I(F_j)=\alpha_j,\quad\deg'_\I K_j^{\pm1}=\deg'_\I K_j'^{\pm1}=\pm\alpha_j,\qquad j\in\I\label{eq:ZI-grade-2}.
    \end{align}

    Let $\tU^+$ be the subalgebra of $\tU$ generated by $E_i$ $(i\in \I)$, let $\tU^-$ be the subalgebra generated by $F_i$ $(i\in\I)$, and let $\tU^0$ be the subalgebra generated by $K_i^{\pm 1},K_i'^{\pm 1}$ $(i\in \I)$. The triangular decomposition is
	\begin{align*} 
		\tU =\tU^+\otimes \tU^0\otimes\tU^-.
	\end{align*}
    For any $\mu=\sum_{i\in\I}m_i\alpha_i\in\Z\I$, we denote $K_\mu=\prod_{i\in\I} K_i^{m_i}$, $K_\mu'=\prod_{i\in\I} (K_i')^{m_i}$.

    The algebra $\tU$ is a Hopf algebra, with the coproduct $\Delta$, the counit $\varepsilon$, and the antipode $S$ defined by
	\begin{align}\label{eq:Hopf structure of tU}
		\begin{split}
			\Delta(E_i)=E_i\otimes 1+K_i\otimes E_i,\quad &\Delta(F_i)=1\otimes F_i+F_i\otimes K_i',
			\\
			\Delta(K_i)=K_i\otimes K_i,\quad &\Delta(K_i')=K_i'\otimes K_i';
			\\
			\varepsilon(E_i)=0=\varepsilon(F_i),\qquad &\varepsilon(K_i)=1=\varepsilon(K_i');
			\\
			S(E_i)=-K_i^{-1}E_i,\quad S(F_i)=-F_iK_i'^{-1},\quad &S(K_i)=K_i^{-1},\quad S(K_i')=K_i'^{-1}.
		\end{split}
	\end{align}

	\begin{lemma}[\cite{Lus93}] 
		\label{lem:anti-involut-QG}
		\quad
		\begin{enumerate}
			\item There exists an involution (called the bar-involution) $u\mapsto \ov{u}$ on $\tU$ given by $\ov{q}=q^{-1}$, $\ov{E_i}=E_i$, $\ov{F_i}=F_i$, and $\ov{K_i}=K_i'$, for $i\in\I$.
			\item 
			There exists an anti-involution $\sigma$ on $\tU$ given by $\sigma(E_i)=E_i$, $\sigma(F_i)=F_i$, and $\sigma(K_i)=K_i'$, for $i\in\I$.
			\item 
			There exists a Chevalley involution $\omega$ on $\tU$ given by $\omega(E_i)=F_i$, $\omega(F_i)=E_i$, and $\omega(K_i)=K_i'$, for $i\in\I$.
		\end{enumerate}
	\end{lemma}

    \begin{lemma} \label{lem:twisting Psi}
		Let $\K^\times=\K\setminus\{0\}$. For scalars $\ba=(a_i)_{i\in\I}\in (\K^{\times})^\I$, there is an automorphism $\widetilde{\Psi}_{\ba}$ of the $\K$-algebra $\tU$ determined by
		\[
		\widetilde{\Psi}_{\ba}:K_i\mapsto a_i^{\frac{1}{2}}K_i,\quad K'_i\mapsto a_i^{\frac{1}{2}}K'_i,\quad E_i\mapsto a_i^{\frac{1}{2}}E_i,\quad F_i\mapsto F_i.
		\]
	\end{lemma}

Recall the Lusztig--Kashiwara skew derivations $\partial_i^L,\partial_i^R: \tU^\pm \rightarrow \tU^\pm$ for $i\in \I$, characterized by
\begin{align}
  E_i y - y E_i &= \frac{K_i \partial_i^L(y)-\partial_i^R(y)K_i'}{q_i-q_i^{-1}},\label{eq:EyyE}\\
  F_i x - x F_i &= \frac{K_i' \partial_i^L(x)-\partial_i^R(x)K_i}{q_i-q_i^{-1}}\label{eq:FxxF}
\end{align}
for $x\in \tU^+$ and $y\in \tU^-$; see \cite[Proposition 3.1.6]{Lus93}. The maps $\partial_i^{L,R}:\tU^-\rightarrow \tU^-$ and $\partial_i^{L,R}:\tU^+\rightarrow \tU^+$ are related by the Chevalley involution $\omega:\tU^\pm\rightarrow \tU^\mp$, that is
\begin{align*}
  \omega(\partial_i^{L,R}(x))=\partial_i^{L,R}(\omega(x)) \qquad \mbox{for all $x\in \tU^+$}.
\end{align*}

\subsection{Braid group symmetries on \texorpdfstring{$\tU$}{U}}
Lusztig introduced braid group symmetries $T'_{i,e},T''_{i,e}$, for $i\in \I$ and $e =\pm 1$, on a quantum group $\mathbf U$ \cite[\S 37.1.3]{Lus93}.
Analogous braid group symmetries $\tTD'_{i,e},\tTD''_{i,e}$, for $i\in \I$ and $e =\pm 1$, exist on the Drinfeld double $\tU$; see \cite[Propositions 6.20--6.21]{LW21b}. The formula for $\tTD'_{i,-1}$ is recalled below.

\begin{proposition}
  [{\cite[Proposition 6.21]{LW21b}}]
   \label{prop:braid1}
 There exists an automorphism $\tTD'_{i,-1}$, for $i\in \I$, on $\tU$ such that
\begin{align*}
\tTD'_{i,-1}(K_{\mu})=&K_{s_i\mu},\qquad \tTD'_{i,-1}(K'_{\mu})=K'_{s_i\mu}, \\
        \tTD'_{i,-1}(E_i)=&-K_i^{-1}F_i,\qquad \tTD'_{i,-1}(F_i)=-E_i(K_i')^{-1},\\
		\tTD'_{i,-1}(E_{j})=&\sum_{r+s=-a_{ij}}(-1)^{r}q_i^{-r}E_i^{(r)}E_{j}E_i^{(s)},\qquad j\neq i,
			\\
		\tTD'_{i,-1}(F_{j})=&\sum_{r+s=-a_{ij}}(-1)^{r}q_i^{r}F_i^{(s)}F_{j}F_i^{(r)},\qquad j\neq i.
			%\text{ for }(j,l)\in\I^{\infty}\text{ and }i\neq(j,l).	
\end{align*}
Moreover, the $\tTD'_{i,-1}$, for $i\in \I$, satisfy the braid relations of the Weyl group $W$.
\end{proposition}

The symmetries $\tTD'_{i,e}$ and $\tTD''_{i,e}$, for $i\in \I$, satisfy the following identities in $\tU$ \cite{LW21b} (analogous to  \cite[37.2.4]{Lus93} in $\mathbf U$)
\begin{align}
  \label{eq:sTs}
  \begin{split}
  \tTD'_{i,-1} &=\sigma \circ \tTD''_{i,+1} \circ \sigma,
  \\
  \tTD''_{i,-e} = \obar \circ \tTD''_{i,+e} \circ \obar,
  & \qquad
  \tTD'_{i,+e} = \obar \circ \tTD'_{i,-e} \circ \obar.
  \end{split}
\end{align}

We sometimes use the following shorthand notation.
\begin{align*}
    \tTD_{i}:= \tTD''_{i,+1},\quad \tTD_{i}^{-1}:= \tTD'_{i,-1}.
\end{align*}
Hence, we define
\begin{align*}  
\tTD_w \equiv \tTD_{w,+1}'' :=\tTD_{i_1}\cdots \tTD_{i_r} \in \Aut (\tU),
\end{align*}
where $w=s_{i_1} \cdots s_{i_r}$ is any reduced expression of $w\in W$.
%Similarly, one defines $T_{w}$ for $w\in W$.

\subsection{${\mathrm i}$Quantum groups}
Given a subset $\bI\subset \I$ of finite type, denote by $W_{\bullet}$ the parabolic subgroup of $W$ generated by $s_i,i\in \bI.$ Set $\bw$ to be the longest element of $W_\bullet$. Let $\Phi_{\bullet}$ be the finite root system generated by simple roots $\alpha_i,i\in \bI.$ Similarly, $\Phi_\bullet^\vee$ is the finite coroot system generated by $\alpha_i^\vee,i\in \bI.$
Let $\rho_\bullet$ be the half sum of positive roots in the root system $\Phi_{\bullet}$, and $\rho_\bullet^\vee$ be the half sum of positive coroots in $\Phi_\bullet^\vee$.

Recall from \cite[Definition 2.3]{Ko14} that an admissible pair $(\I=\wI\cup \bI,\tau)$ of Kac--Moody type consists of a partition $\wI\cup \bI$ of $\I$ and a Dynkin diagram involution $\tau$ of $\fg$ such that
\begin{itemize}
\item[(1)] $\bI$ is of finite type,

\item[(2)]
 $\bw(\alpha_j) = - \alpha_{\tau j}$ for $j\in \bI$,
\item[(3)]
 If $j\in \wI$ and $\tau j =j$, then $\alpha_j(\rho_{\bullet}^\vee)\in \Z$.
\end{itemize}

Given a symmetric pair $(\I=\wI\cup\bI,\tau )$, the universal iquantum group $\tUi$ is the $\K$-subalgebra of $\tU$ generated by
	\begin{align}
\label{def:iQG}
\begin{split}
&B_i=F_i+\tTD_{\bw}( E_{\tau i}) K_i', \qquad \tk_i^{\pm1}=(K_i K'_{\tau i})^{\pm1} \qquad (i\in \wI),
\\
&E_j,F_j,K_j^{\pm1},K'^{\pm1}_j\qquad (j\in \bI).
\end{split}
\end{align}
 Let $\tU^{\imath 0}$ be the $\K$-subalgebra of $\tUi$ generated by $\tk_i^{\pm1},K_j^{\pm1},K'^{\pm1}_j$, for $i\in \wI,j\in\bI$. Let $\tbU$ denote the quantum group associated with $\bI$; it is a subalgebra of $\tUi$.
The algebra $\tUi$ is a right coideal subalgebra of $\tU$; the pair $(\tU,\tUi)$ is called a quantum symmetric pair, and $\tUi$ is called the universal iquantum group.
\begin{lemma}[{\cite[Proposition 3.12]{WZ23}}]
There exists an anti-involution $\sigma^\imath$ on $\tUi$ given by $\sigma^\imath(B_i)=B_i$, $\sigma^\imath(x)=\sigma(x)$, for $i\in\wI,x\in\tU^{\imath 0}\tbU$.  
\end{lemma}

\begin{lemma}[cf. {\cite[Lemma 3.4]{Ko14}}]\label{lem:Tw-on-X} 
 For all $i\in \I_\bullet$, one has
  \begin{align*}
    \tTD_{w_\bullet}(E_i)&=-F_{\tau i} K_{\tau i}'^{-1},\qquad
    \tTD_{w_\bullet}(F_i)=-K_{\tau i}^{-1}E_{\tau i},\\
    \tTD_{w_\bullet}^{-1}(E_i)&=-K_{\tau i}^{-1} F_{\tau i},\qquad
    \tTD_{w_\bullet}^{-1}(F_i)=-E_{\tau i} K_{\tau i}'^{-1}.
  \end{align*}
\end{lemma}
%%%%%%%%%%%%%
\subsection{Letzter map}\label{sec:letzter map}
Define the left and right adjoint actions $\ad_l$ and $\ad_r$ of $\tU$ by
\begin{equation}\label{adjoint action}
    \ad_l(a)(b)=a_{(1)}bS(a_{(2)}),\quad \ad_r(a)(b)=S(a_{(1)})ba_{(2)},
\end{equation}
where $\Delta(a)=a_{(1)}\otimes a_{(2)}$ is written in Sweedler notation.

Set $H:=\tbU\tU^0$ and define two subalgebras of $\tU$ by
\begin{align}
    \cR_{\bullet}^{-,r}:=\langle\ad_r(H)(F_i)\mid i\in\wI\rangle,\qquad\cR_{\bullet}^{+,r}:=\langle \ad_r(H)(E_iK_i^{-1})\mid i\in\wI\rangle.
\end{align}
The arguments of \cite[\S 2.2]{KY26} give the following linear isomorphisms, induced by multiplication:
\begin{align}
    \label{eq:triang_U of R}
     \cR_\bullet^{-,r}\otimes H &\otimes \cR_\bullet^{+,r}\cong \tU \\
     \label{eq:triang_Upm of R}
    \cR_\bullet^{-,r}\otimes \tbU^-\cong \tU^-, &\qquad S(\tbU^+)\otimes \cR_\bullet^{+,r} \cong S(\tU^+).
\end{align}
It follows that $\cR_\bullet^{\pm,r}$ are right coideal subalgebras of $\tU$. The subalgebra $\cR_\bullet^{-,r}$ can also be described in terms of braid group
automorphisms or skew derivations.

\begin{proposition}[{\cite[Proposition 2.4, Lemma 2.5]{KY26}, \cite[Lemma 2.4]{KY21}}]\label{another description of R}
  The following descriptions hold:
  \begin{align}
      \cR_\bullet^{-,r}=&\bigcap_{j\in \I_\bullet} \big(\tU^-\cap \tTD_j(\tU^-)\big)=\{u\in \tU^-\,|\,\partial_j^R(u)=0\quad \mbox{for all $j\in \bI$}\}\label{R of derivation}\\
      =&\tU^-\cap \tTD_{\bw}(\tU^-) \label{R of T}
  \end{align}
\end{proposition}
\begin{remark}
 The descriptions in Proposition \ref{another description of R} only require $\bI$ to be of finite type.   
\end{remark}

Consider the subalgebra $\cA\subset \tU$ defined by
\begin{align*}
  \cA=\K\langle F_i, E_j, \tU^{\imath 0}\mid i\in \I, j\in \bI\rangle.
\end{align*}
$\cA$ is called the horospherical subalgebra corresponding to $(\I=\wI\cup\bI,\tau)$. Also set $\cH_\bullet:=\tbU\tU^{\imath 0}\subset \tU$.
By \eqref{eq:triang_Upm of R} and \cite[Theorem 2.1]{KY21}, the multiplication maps
\begin{align}\label{eq:A-triang-1}
    \cR_\bullet^{-,r}\otimes \cH_\bullet\longrightarrow \cA,\qquad \cH_\bullet\otimes \cR_\bullet^{-,r} \longrightarrow \cA
\end{align}
are isomorphisms. 

Consider the subalgebra $\tU^{\text{poly}}\subset \tU$ defined by
\begin{equation*}
    \tU^{\text{poly}}=\langle\cA,E_jK_j^{-1},K_j^{-1}\mid j\in\wI\rangle.
\end{equation*}
The arguments of \cite[\S 2.5]{KY26} yield
\begin{align}
    \label{eq:Upoly-triang-1}
    \cA\otimes \K[K_i^{-1}\,|\,i\in \I_\circ]\otimes \cR_{\bullet}^{+,r}\cong \tU^\poly.
\end{align}
Set $\cR_{\bullet,+}^{+,r}=\cR_\bullet^{+,r}\cap\ker\eps$ and consider the left ideal $\cJ_\bullet$ of $\tU^\poly$ defined by
\begin{equation}\label{eq:J-triang}
   \cJ_\bullet=\tU^\poly\Bigl(\cR_{\bullet,+}^{+,r}
+\sum_{j\in\wI}\K K_j^{-1}\Bigr)= \tU^\poly\Span_\K\{E_jK_j^{-1},K_j^{-1}\mid j\in\wI\}\cH_\bullet.
\end{equation}
The triangular decomposition \eqref{eq:Upoly-triang-1} yields the vector space decomposition
  \begin{align}\label{eq:Upoly-decomp}
     \tU^\poly = \cA\oplus \cJ_{\bullet}.
  \end{align}

Let $\psi:\tU^\poly\rightarrow \cA$ denote the $\K$-linear projection with respect to the direct sum decomposition \eqref{eq:Upoly-decomp}. Since $\tUi$ is a subalgebra of $\tU^{\text{poly}}$, $\psi$ restricts to a map $\psi:\tUi\rightarrow \cA$, called the \textit{Letzter map}.

The algebra $\tUi$ has a natural $\N$-filtration $\cF_\ast$ defined via the degree function given by 
\begin{equation}
\begin{aligned}\label{eq:filt-def}
  \deg(B_i)=1 \quad \mbox{for $i\in \wI$},\qquad
  \deg(h)=0 \quad \mbox{for $h\in \cH_\bullet$.}
\end{aligned}
\end{equation}
The horospherical subalgebra $\cA$ is $\N$-graded via the degree function given by
  \begin{align*}
     \deg(F_i)=1 \quad \mbox{for all $i\in \wI$,}\qquad
     \deg(h)=0 \quad \mbox{for all $h\in \cH_\bullet$}.
  \end{align*}
For $n\in \N$, let $\cA_n$ denote the $n$-th graded component of $\cA$ and write $\cA_{<n}=\bigoplus_{m=0}^{n-1}\cA_m$.
Let $\gr(\tUi)$ be the associated graded algebra. Following \cite[Lemma 2.7]{KY26}, there exists an algebra isomorphism
\begin{align}\label{eq:varphi}
  \varphi:\cA \rightarrow \gr(\tUi)
\end{align}
given by $\varphi(F_i)=B_{i}$ for all $i\in \wI$ and $\varphi(h)=h$ for all $h\in \cH_\bullet$.

By construction, the Letzter map $\psi:\tUi\rightarrow \cA$ is a linear map of filtered vector spaces. Let
\begin{align*}
  \gr(\psi):\gr(\tUi)\rightarrow \gr(\cA)\cong \cA
\end{align*}
be the associated graded map. Its composition with $\varphi$ satisfies $\gr(\psi)\circ \varphi=\id_\cA$.

\begin{lemma} [{\cite[Lemmas 2.10 and 2.11]{KY21}}]
    \label{lem:Letzter-props}
    The Letzter map $\psi:\tUi\rightarrow \cA$ is a linear isomorphism of filtered vector spaces with associated graded map $\gr(\psi)=\varphi^{-1}$. Moreover, the following hold:
    \begin{enumerate}
      \item $\psi(ab) = \psi(a\psi(b))$ for all $a,b\in \tUi$;
      \item $\psi(h_1bh_2) = h_1\psi(b) h_2$ for all $b\in \tUi$, $h_1,h_2\in \cH_\bullet$;
      \item $\psi(\tTD_i(b))=\tTD_i(\psi(b))$ for all $i\in \bI$;
      \item $\psi(ab)-\psi(a)\psi(b)\in \cA_{<m+n}$ for all $a\in \cF_m(\tUi)$, $b\in \cF_n(\tUi)$.  
    \end{enumerate}
\end{lemma}
Note that Lemma \ref{lem:Letzter-props}~(1)-(3) are valid for $b\in\tU^\poly$.
\iffalse
\begin{remark}
If we choose $\wI'\subseteq\wI$ such that $\tau(\wI')=\wI'$, then for $(\I'=\wI'\cup\bI,\tau|_{\I'})$, we can define iquantum group $\tUi_s$ and $\psi_s:\tU^\poly_s\rightarrow\cA_s$. By \eqref{eq:Upoly-decomp}, we can find $\psi_s=\psi|_{\tU^\poly_s}$. Since Letzter map is isomorphism, we can know that $\psi^{-1}|_{\tUi_s}=\psi_s^{-1}$.
\end{remark}
\fi
\begin{remark}\label{rem:local Letzter and Letzter}
Let $\wI'\subseteq\wI$ be $\tau$-stable, and set
$\I'=\wI'\cup\bI$. Denote by $\tUi_s$ and $\cA_s$ the
iquantum group and the horospherical algebra associated with
the Satake subdiagram $(\I',\tau|_{\I'})$, respectively, and let $\psi_s:\tU_s^\poly\longrightarrow\cA_s$
be the corresponding projection. By \eqref{eq:Upoly-decomp}, $\psi_s=\psi|_{\tU_s^\poly}$.
Since the Letzter maps $\psi:\tUi\rightarrow\cA$ and $\psi_s:\tUi_s\rightarrow\cA_s$ are linear isomorphisms, we obtain 
$\psi^{-1}|_{\cA_s}=\psi_s^{-1}$.
\end{remark}
\subsection{Star product of \texorpdfstring{$\cA$}{A}}
Assume that $A$ is an $\N$-graded $\K$-algebra $A=\bigoplus_{j\in \N} A_j$. For all $n, m \in \N$, we denote 
\[
A_{[n,m]} : = 
\begin{cases}
A_n \oplus \ldots \oplus A_m &\mbox{if} \; \; n \leq m,
\\
0 &\mbox{otherwise},
\end{cases}
\]
and $A_{<m} := A_{[0, m-1]}$, $A_{\le m}:= A_{[0,m]}$. 

\begin{definition}
[{\cite[Subsection 2.1]{ES20}}]\label{def:star-prod}
Let $A$ be an $\N$-graded $\K$-algebra.
\begin{enumerate}
\item[(i)] A star product on $A$ is an associative bilinear operation $* : A \times A \to A$, $(a,b)\mapsto a\ast b$
such that
\begin{align}\label{eq:a*b-ab}
  a * b - ab \in A_{< m+ n} \qquad \mbox{for all $a \in A_m, b \in A_n$.}
\end{align}
\item[(ii)] A star product $*$ on $A$ is called short if 
\[
a * b \in A_{[|m-n|, m+n]} \qquad \mbox{for all $a \in A_m, b\in A_n$.} 
\]
\item[(iii)] A star product $\ast$ on $A$ is called $0$-equivariant
if
\[
a * h = ah \quad \mbox{and} \quad  h * a= h a \qquad \mbox{for all $h \in A_0, a \in A$.} 
\]
\item[(iv)] A star product $\ast$ on $A$ is called $\Z/2$-equivariant if 
\begin{equation}
\label{Z2cond}
a * b \in \bigoplus_{k=0}^{\lfloor \frac{m+n}{2} \rfloor} A_{m+n-2 k} \qquad \mbox{for all $a \in A_m, b \in A_n$.} 
\end{equation}
\end{enumerate} 
\end{definition}

The arguments of \cite[\S\S~3.1--3.2]{KY26} show that
\begin{align}\label{eq:ast-def}
  a\ast b = \psi(\psi^{-1}(a)\psi^{-1}(b)) \qquad \text{for all $a,b\in \cA$}
\end{align}
defines a star product on $\cA$. By Lemma \ref{lem:Letzter-props}, the star product $\ast$ on $\cA$ is $0$-equivariant. Moreover, \cite[Lemma 3.7, Theorem 3.10]{KY26} shows that the star product $\ast$ on $\cA$ is $\Z/2$-equivariant and short.

%\begin{proposition}[{\cite[Lemma 3.7, Theorem 3.10]{KY26}}]
%The star product $\ast$ on $\cA$ is $\Z/2$-equivariant and short.    
%\end{proposition}

%\red{
\begin{proposition}[{\cite[Theorem 4.6]{KY26}}]\label{sigma of (A,*)}
The anti-involution $\sigma$ of $\tU$ restricts to an  anti-automorphism of $(\cA,\ast)$.    
\end{proposition}
%}
%%%%%%%%%%%%%%%%
\section{\texorpdfstring{$\mathrm i$}{i}Braid group symmetries}\label{sec:i braid action}
In this section, we prove \cite[Conjecture 5.13]{WZ23} and establish the ibraid group symmetries on $\tUi$.
\subsection{Relative Weyl groups}
Given an admissible pair $(\I=\bI \cup \wI,\tau)$, the corresponding involution $\theta$ (acting on the weight lattice) is recovered as
\begin{align}
  \label{eq:theta}
\theta=-\bw \circ \tau.
\end{align}
Fix a set $\I_{\circ,\tau}$ of representatives of the $\tau$-orbits in $\wI$. For $i\in \I_{\circ,\tau}$, set $\I_{\bullet,i}:=\bI\cup \{i,\tau i\}$ and there is a rank one Satake subdiagram $(\I_{\bullet,i}=\{i,\tau i\}\cup\I_\bullet,\tau|_{\I_{\bullet,i}})$. Let $\fwItau$ be the subset of $\wItau$ given by
\begin{align}
\label{def:fwItau}
\fwItau=\{i\in \wItau| \; \I_{\bullet,i} \text{ is of finite type}\}.
\end{align}

  For $i\in\fwItau$, let $W_{\bullet,i}$ be the parabolic subgroup of $W$ generated by $s_j,j\in \I_{\bullet,i}$ and let $w_{\bullet,i}$ be the longest element of $W_{\bullet,i}$. Denote by $\ell$ the length function of the Coxeter system ($W_{\bullet,i},\bIi$).
Define $\br_i\in W_{\bullet,i}$ such that
\begin{align}\label{def:bsi}
\bwi= \br_i w_\bullet \, (=w_\bullet \br_i),
\qquad \text{ where }
\ell(\bwi) = \ell(\br_i) + \ell(w_\bullet).
\end{align}
Denote by $\tau_{\bullet,i}$ the diagram involution on $\bIi$ such that
\begin{align} \label{def:taui}
\bwi(\alpha_j) =-\alpha_{\tau_{\bullet,i} j}, \qquad \forall j\in \bIi.
\end{align}
The subgroup of $W$,
\[
W^\circ :=\langle \br_i\mid i\in \fwItau \rangle,
\]
is itself a Weyl group with simple reflections identified with $\{\br_i \mid i\in \fwItau\}$.

The {\em relative Weyl group} associated to the Satake diagram $(\I=\bI \cup \wI,\tau)$ can be identified with $W^\circ$. Let $\{\ov\alpha_i\mid i\in \fwItau\}$  be the simple system of the relative (or restricted) root system, where $\ov\alpha_i$ is identified with the following element (cf. \cite[\S 2.3]{DK19})
\begin{align}
 \label{def:balpha}
\ov\alpha_i:=\frac{\alpha_i - \theta(\alpha_i)}{2},\qquad (i\in \fwItau).
\end{align}
Note that $\ov\alpha_i=\ov\alpha_{\tau i}$.

We shall refer to the braid group associated to the relative Weyl group $W^\circ$ as the {\em relative braid group} and denote it by $\mathrm{Br}(W^\circ)$. Accordingly, we denote the braid group associated to $W_\bullet$ by $\mathrm{Br}(W_\bullet)$.  It is well known that (see, e.g., \cite[\S 2.2]{DK19})
\begin{equation}\label{eq:small weyl group}
    W_\bullet \rtimes W^\circ\cong W^\theta=\{w\in W\mid w\circ\theta=\theta\circ w\}.
\end{equation}

\subsection{Quasi $K$-matrix}

%The relative braid (or ibraid) group symmetries $\tTT_{i,e}'$ and $\tTT''_{i,e}$ ($i\in\I$, $e\in\{+1,-1\}$) on $\tUi$ are established in \cite{WZ23}.

Viewing $B_i$ as an element of $\tU$, set 
$$B_i^\sigma=\sigma(B_i)=F_i+K_i\tTD^{-1}_{\bw}( E_{\tau i})\text{ for }i\in\I_\circ.$$ 
\begin{proposition}
		[\text{\cite[Theorem~3.6]{WZ23}}]
		\label{prop:Kmatrix}
There exists a unique element $\tfX=\sum_{\mu \in \N \I} \tfX^{\mu}$ such that $\tfX^0=1, \tfX^{\mu}\in \tU_{\mu}^+$ and the following intertwining relations hold:
 \begin{align}\label{eq:fX2}
 \begin{split}
 B_i  \tfX &=  \tfX B_i^{\sigma},\qquad (i\in \wI),\\
 x  \tfX &= \tfX x,\qquad (x\in \tU^{\imath 0}\tbU).
 \end{split}
 \end{align}
 Moreover, $\tfX^\mu =0$ unless $\theta (\mu) =-\mu$.
 \end{proposition}

The element $\widetilde{\Upsilon}$ is called the quasi $K$-matrix; see \cite{BW18b,BK19,AV22}.
Denote by $\tU_{\bullet,i}$ the quantum group associated to $ \I_{\bullet,i}$. Let $\widetilde{\Upsilon}_i$ be the rank-one quasi $K$-matrix associated to $ \I_{\bullet,i}$, i.e., $\widetilde{\Upsilon}_i=\sum_{m\in \N}\widetilde{\Upsilon}_{i,m}$ with $\widetilde{\Upsilon}_{i,m}\in\tU_{m(\alpha_i-\theta\alpha_i)}^+$, and $\widetilde{\Upsilon}_{i}^{0}=1$ by \cite[Remark 3.10]{WZ23}. Define a distinguished parameter $\bvs_\diamond=(\varsigma_{i,\diamond})_{i\in\I}$ by
\begin{align}
		\label{eq:disting-para}
		\varsigma_{i,\diamond}= \left\{
		\begin{array}{ll}
			1, & \text{ if } i\in\bI,
			\\
			- q^{-(\alpha_i , \alpha_i+\bw\alpha_{\tau i})/2}, & \text{ if } i\in\wI.
		\end{array}
		\right.
	\end{align}
    Recall the automorphism $\widetilde{\Psi}_{\varsigma_{\diamond}}$ of $\tU$ from Lemma~\ref{lem:twisting Psi}, and set
	\begin{equation}\label{eq: T_i twisted def}
		\tT_i:=\widetilde{\Psi}_{\bvs_\diamond}^{-1}\circ\tTD_i\circ \widetilde{\Psi}_{\bvs_\diamond},
		\qquad \tT_i^{-1}:=\widetilde{\Psi}_{\bvs_\diamond}^{-1}\circ\tTD_i^{-1}\circ \widetilde{\Psi}_{\bvs_\diamond}.
	\end{equation}
	Clearly $\tT_i$ and $\tT_i^{-1}$, for $i\in\I$, are automorphisms of $\tU$ and satisfy the braid group relations. Hence, we define
	$\tT_w:=\tT_{i_1}\cdots \tT_{i_r},$
	where $w=s_{i_1}\cdots s_{i_r}$ is any reduced expression.

    For convenience, we set $\tk_i^{\pm1}=(K_i K'_{\tau i})^{\pm1}$ where $i\in\bI$.
For $\lambda=\sum_{i\in \I} m_i \alpha_i \in \Z\I$, we use the notation
\begin{align}
  \label{eq:vkHla}
\vs_{\lambda,\dm}:=\prod_{i\in \I}\vs_{i,\dm}^{m_i},\qquad
\tk_{\lambda}:=\prod_{i\in \I}\tk_{i}^{m_i}  \in \tU^{\imath 0},\qquad
\cK_\lambda=K_\lambda K_{-\theta\lambda}'.
\end{align}

\begin{proposition}[\text{cf. \cite[Theorem 4.2]{BW18b}}]
  \label{prop:Tjblack}
Let $i\in \bI$. The automorphism $\tT'_{i,-1} =\tTD'_{i,-1}$ on $\tU$ restricts to an automorphism of $\tUi$. Moreover, the action of $\tT'_{i,-1}$ on $B_j \; (j\in \wI)$ is given by
\begin{align}\label{eq:newb-1}
\begin{split}
\tT'_{i,-1}(B_j)=\sum_{s=0}^r (-1)^s q_i^s F_i^{(r-s)} B_j F_i^{(s)},\qquad \text{ for } r=-a_{ij}.
\end{split}
\end{align}
\end{proposition}

\subsection{Characterization of $\tT_i^{-1}(B_j)$}
Recall the right adjoint actions $\ad_r$ of $\tU$ in \eqref{adjoint action}. Let $\cM_{\bullet,i,j}:=\ad_r(\tbiU^-)(F_j)$ be the subspace for $i\in\fwItau$ and $j\notin \bIi$. Let $\tbiU^\imath$ be the iquantum group  associated with $(\bIi=\{i,\tau i\}\cup\bI,\tau|_{\bIi})$ which is the subalgebra of $\tUi$. Set
 $$\tU^\imath_{\bullet,i,+}=\langle\tbU,B_{j},\tk_j\mid j\in\{i,\tau i\}\rangle.$$
 %\red{$\cK_j$ not invertible?}
 Recall $\cK_\lambda=K_\lambda K_{-\theta\lambda}'$ in \eqref{eq:vkHla}. Note that $\cK_j=\tk_jK'_{\bw\alpha_{\tau j}-\alpha_{\tau j}}\in\tU^{\imath 0}$.
\begin{proposition}\label{prop:Mij-subset-BFB}
  For $i\in\fwItau$ and $j\notin \bIi$, we have $\cM_{\bullet,i,j}\subseteq\tU^\imath_{\bullet,i,+}F_j\tU^\imath_{\bullet,i,+}$.
  %and $\sigma(\cM_{\bullet,i,j})\subseteq\sigma(\tU^\imath_{\bullet,i,+})F_j\sigma(\tU^\imath_{\bullet,i,+})$.  
\end{proposition}
\begin{proof}
By $\ad_r(E_s)(F_j)=0$ for $s\in\bIi$, we know that $\cM_{\bullet,i,j}$ is a $\tbiU$-module via $\ad_r$. Since $\tU$ is a Hopf algebra, for $x,b\in\tU$, we have 
   \begin{equation}\label{eq:Hopf identity}
xb=\varepsilon(b_{(1)})xb_{(2)}=b_{(1)}S(b_{(2)})xb_{(3)}=b_{(1)}\ad_r(b_{(2)})(x).
   \end{equation}
By \cite[Lemma 2.6]{KY26}, 
\begin{equation*}
 \Delta(B_i)=B_i\otimes K_i'+1\otimes F_i+\sum_{0\ne\nu\in\N\bI}h_\nu K_{-\nu}\cK_{-\theta\alpha_i}\otimes u_{-\theta\alpha_i-\nu}K_i',
\end{equation*}
where $h_\nu\in\tbU^+,u_{-\theta\alpha_i-\nu}\in(\tbiU^+)_{-\theta\alpha_i-\nu }$. Using \eqref{eq:Hopf identity}, we get
\begin{equation}\label{eq:apply Hopf identity to xBi}
  \ad_r(F_i)(x)=xB_i-B_i\ad_r(K_i')(x)-\sum_{0\ne\nu\in\N\bI}h_\nu K_{-\nu}\cK_{-\theta\alpha_i}\ad_r( u_{-\theta\alpha_i-\nu}K_i')(x). 
\end{equation}

We verify $x\in \tU^\imath_{\bullet,i,+}F_j\tU^\imath_{\bullet,i,+}$ for $x\in\cM_{\bullet,i,j}$ by induction on the height of $\wt(x)$. For $x=F_j$, it is obvious. Assume now $x\in\tU^\imath_{\bullet,i,+}F_j\tU^\imath_{\bullet,i,+}$ for $x\in\cM_{\bullet,i,j}$. It suffices to prove $\ad_r(F_s)(x)\in\tU^\imath_{\bullet,i,+}F_j\tU^\imath_{\bullet,i,+}$ for $s\in\bIi$. When $s\in\bI$, it is obvious by $\tbU\subseteq\tUi$. When $s\in\{i,\tau i\}$, using \eqref{eq:apply Hopf identity to xBi}, we can find that $\ad_r(F_s)(x)\in\tU^\imath_{\bullet,i,+}F_j\tU^\imath_{\bullet,i,+}$.
\iffalse
 \red{  So for $x\in\cM_{\bullet,i,j}$, we prove $x\in \tU^\imath_{\bullet,i,+}F_j\tU^\imath_{\bullet,i,+}$ by induction on the height of $\wt(x)$. It is obvious for $x=F_j$. For $\ad_r(F_l)(x)$ where $l\in\bIi$, it is obvious for $l\in\bI$. For $l\in\{i,\tau i\}$, using \eqref{eq:apply Hopf identity to xBi}, we complete the induction.}
 \fi
\end{proof}

\begin{corollary}
     \label{cor:rktwoRPijb}
Let $i\in\fwItau$ and $j\notin \bIi$. Then there exist non-commutative polynomials $R_{ij,\bullet}(x_i,x_{\tau i}, y_i, y_{\tau i},z)$ and $P_{ij,\bullet}(x_i, x_{\tau i}, y_i, y_{\tau i}, z)$ with coefficients in $\tbU$, each of which is linear in $z$, such that
\begin{align}
  \tT_{\br_i}^{-1} (F_j) =& R_{ij,\bullet}(B_i^\sigma, B_{\tau i}^\sigma, \cK_i, \cK_{\tau i}, F_j),\label{eq:TF-poly-1}\\
  \tT_{\br_i}(F_j) =& P_{ij,\bullet}(B_i, B_{\tau i}, \cK_i, \cK_{\tau i}, F_j)\label{eq:TF-poly-2}.
\end{align}
\end{corollary}
\begin{proof}
    Since $\sigma(\tT_{\br_i}(F_j))=\tT_{\br_i}^{-1}(F_j)$, we only need to prove \eqref{eq:TF-poly-2}. By $\br_i\alpha_j>0$, we know $\tT_{\br_i}(F_j)=\tTD_{\br_i}(F_j)\in\tU^-$. Then $\tTD_{\br_i}(F_j)\in\cM_{\bullet,i,j}$ by \cite[Lemma 39.1.2]{Lus93}. The desired formula \eqref{eq:TF-poly-2} follows from Proposition \ref{prop:Mij-subset-BFB}.
\end{proof}

Recall the antipode $S$ in \eqref{eq:Hopf structure of tU}. We define the following algebra automorphisms
\begin{equation}\label{eq:auto-Omega}
 \Ob:=\tTD_{\bw}\circ\tau\circ S\circ\sigma\circ\omega,\qquad \wObb:=\tPsi_{\bvs_{\dm}}^{-1}\circ\Omega_{\bullet}\circ\tPsi_{\bvs_{\dm}}.
\end{equation}
In fact, using Lemma~\ref{lem:Tw-on-X} and \eqref{eq:Hopf structure of tU},
a direct check on generators implies
\begin{equation}\label{eq:S-sigma-omega=Tbw-tau}
S\circ\sigma\circ\omega|_{\tbU}=\tTD_{\bw}\circ\tau|_{\tbU}.   
\end{equation}
Using \eqref{eq:Hopf structure of tU} and Lemma~\ref{lem:twisting Psi}, for $i\in\bI,j\in\wI$, we have
\begin{equation}\label{eq:formula of wObb}
 \begin{aligned}
 \wObb(K_i)=K_i,\quad  \wObb(K_i')=K_i' \quad \wObb(E_i)&=q_i^{-2}E_i\cK_i^{-1},\quad \wObb(F_i)=q_i^{-2}F_i\cK_i,\\
 \wObb(F_j)=-q_j^{-2}\tTD_{\bw}( E_{\tau j}) K_j'\cK_{\theta\alpha_j},\qquad&\wObb(K_j\tTD_{\bw}^{-1}( E_{\tau j}))=-q_j^{-2}F_j\cK_{\theta\alpha_j}.
\end{aligned}   
\end{equation}
\begin{lemma}[{\cite[Lemma 7.1]{BK19}}]
  For $i\in\I$, on $\tU$ we have
  \begin{equation}\label{S-sigma-omega commute T}
S\circ\sigma\circ\omega\circ\tTD_i=\tTD_i\circ S\circ\sigma\circ\omega.
  \end{equation}
\end{lemma}
Since $\tau\circ\tTD_j=\tTD_{\tau j}\circ\tau$ for $j\in\I$ and $\tau(\br_i)=\br_i$, using \eqref{S-sigma-omega commute T}, we get
\begin{equation}\label{eq:wObb commute Tbr}
   \wObb\circ \tT_{\br_i}^{\pm1}=\tT_{\br_i}^{\pm1}\circ\wObb.
\end{equation}

Recall the algebra $\tU$ is bigraded by \eqref{eq:ZI-grade-1} and \eqref{eq:ZI-grade-2}. For $\mu,\nu\in\Z\I$, let $\tU_{(\mu,\nu)}$ denote the bihomogeneous subspace $\{x\in\tU\mid \deg_\I(x)=\mu,\ \deg_\I'(x)=\nu\}$, so that $\tU=\oplus_{\mu,\nu\in\Z\I} \tU_{(\mu,\nu)}$. For $\lambda\in\Z\I$, we define
\begin{equation}\label{eq:sum of bigraded subspace and int subalge}
  \tU[\lambda]:=\bigoplus_{\substack{\mu,\nu\in\Z\I,\\2\lambda=2\nu-(1+\theta)\mu}}\tU_{(\mu,\nu)},\qquad \tU_{\text{int}}:=\bigoplus_{\lambda\in\Z\I}\tU[\lambda].
\end{equation}
Note that the bigraded subspace $\tU_{\text{int}}$ is a bigraded subalgebra of $\tU$. Obviously, we know $\tbU\subseteq\tU_{\text{int}}$. A direct check implies $\{F_j,\tTD_{\bw}( E_{\tau j}) K_j',K_j\tTD_{\bw}^{-1}( E_{\tau j}),\cK_j\mid j\in\wI \}\subseteq\tU_{\text{int}}$. For example, for $j\in\wI$, we have
\begin{equation*}
  \deg_\I(\cK_j)=\alpha_j-\theta\alpha_j=\deg_\I'(\cK_j), \quad 2(1-\theta)\alpha_j-(1-\theta)(1+\theta)\alpha_j=2(1-\theta)\alpha_j.
\end{equation*}
Thus, $\cK_j\in\tU[(1-\theta)\alpha_j]$. So we get
\begin{equation}\label{eq:Ui and sigmaUi in Uint}
    \tUi\cup\sigma(\tUi)\subseteq\tU_{\text{int}}.
\end{equation}
\begin{remark}
    We define a lattice $L$ of $\Z\I$ by
    \begin{equation*}
      L=\{\mu\in\Z\I\mid (1+\theta)\mu\in2\Z\I\}=\{\nu\in\Z\I\mid \theta\nu=\nu\}\oplus\{\nu\in\Z\I\mid \theta\nu=-\nu\}.  
    \end{equation*}
    Then we have $\tU_{\text{int}}=\bigoplus_{\nu\in\Z\I,\mu\in L}\tU_{(\mu,\nu)}$. Set
    \begin{equation*}
      \tU_{\text{int}}^0:= \tU_{\text{int}}\cap\tU^0=\Span_{\K}\{K_\alpha K'_\beta\mid\alpha+\beta\in L\}. 
    \end{equation*}
    Since $\tU^-\otimes\tU^0\otimes S(\tU^+)\cong\tU$ by multiplication map, we can find that the multiplication map
    \begin{equation}\label{eq:structure of U-int}
     \tU^-\otimes\tU_{\text{int}}^0\otimes S(\tU^+)\longrightarrow\tU_{\text{int}}   
    \end{equation}
    is an isomorphism. Note that $\tU_{\text{int}}$ is a right coideal subalgebra of $\tU$.
\end{remark}

By Proposition~\ref{prop:braid1} and \eqref{eq:Hopf structure of tU}, for $x\in\tU$ with $\deg_\I(x)=\mu,\ \deg_\I'(x)=\nu$, and $i\in\I$, a direct check on generators implies
\begin{align}
&\deg_\I(\tau\circ S\circ\sigma\circ\omega(x))=-\tau(\mu),\quad&&\deg'_\I(\tau\circ S\circ\sigma\circ\omega(x))=-\tau(\nu),\label{eq:tau-S-sig-omega retain grade}\\
  &\deg_\I\tTD_i^{\pm1}(x)=s_i(\mu),\quad && \deg'_\I\tTD_i^{\pm1}(x)=s_i(\nu).\label{eq:T retain grade}
\end{align}

By \eqref{eq:tau-S-sig-omega retain grade}, \eqref{eq:T retain grade} and \eqref{eq:sum of bigraded subspace and int subalge}, we know that $\wObb$ is an algebra endomorphism of $\tU_{\text{int}}$. Since $\cK_{-\theta\lambda}\in\tU[(1-\theta)\lambda]$, we define a linear map on $\tU_{\text{int}}$ 
\begin{equation}\label{eq:formula of Xi}
  \Xi:\tU_{\text{int}}\longrightarrow\tU_{\text{int}},\qquad \Xi(x)=\vs_{\lambda,\dm}^{-1}q^{\frac{1}{2}(\lambda,\lambda+\theta\lambda)}\wObb(x)\cK_{-\theta\lambda},\ x\in\tU[\lambda].
\end{equation}
 
\begin{lemma}\label{lem:Xi is alge-homomorphism}
   The linear map $\Xi:\tU_{\text{int}}\longrightarrow\tU_{\text{int}}$ is an algebra endomorphism of $\tU_{\text{int}}$. Moreover, we have
   \begin{equation}\label{eq:formula of Xi on generators}
      \Xi|_{\tbU}=\id_{\tbU},\quad \Xi(F_i)=\tTD_{\bw}(E_{\tau i})K_i',\quad \Xi(B_i^\sigma)=B_i,\quad \Xi(\cK_i)=\cK_i,\qquad i\in\wI.
   \end{equation}
\end{lemma}
\begin{proof}
  By \eqref{eq:wt-grade}, \eqref{eq:ZI-grade-1} and \eqref{eq:ZI-grade-2}, for $\mu,\nu\in\Z\I$ and $\deg_\I(x)=\mu,\ \deg_\I'(x)=\nu$, we know  $\wt(x)=\mu-\nu$. Then for $x\in\tU[\lambda]\cap\tU_{(\mu,\nu)}$ and $y\in\tU[\lambda']\cap\tU_{(\mu',\nu')}$, we have
  \begin{align*}
   \Xi(x)\Xi(y)=& \vs_{\lambda+\lambda',\dm}^{-1}q^{\frac{1}{2}(\lambda,\lambda+\theta\lambda)+\frac{1}{2}(\lambda',\lambda'+\theta\lambda')}\wObb(x)\cK_{-\theta\lambda}\cdot \wObb(y)\cK_{-\theta\lambda'}\\
   =&\vs_{\lambda+\lambda',\dm}^{-1}q^{\frac{1}{2}(\lambda,\lambda+\theta\lambda)+\frac{1}{2}(\lambda',\lambda'+\theta\lambda')-(\theta\mu'-\theta\nu',\lambda+\theta\lambda)}\wObb(xy)\cK_{-\theta(\lambda+\lambda')}.
  \end{align*}
  For $\alpha,\beta\in\Z\I$, we know $\big((\theta-1)\alpha,(\theta+1)\beta\big)=0$ since 
  \[\big((\theta-1)\alpha,(\theta+1)\beta\big)=\big(\theta(\theta-1)\alpha,\theta(\theta+1)\beta\big)=-\big((\theta-1)\alpha,(\theta+1)\beta\big).\]
  Since $2\lambda'=2\nu'-(1+\theta)\mu'$, we know\
  \begin{equation*}
   -(\theta\mu'-\theta\nu',\lambda+\theta\lambda)=(\nu'-\theta\mu',\lambda+\theta\lambda)=\big(\lambda'-\frac{1}{2}(\theta-1)\mu',(\theta+1)\lambda\big)=(\lambda',\lambda+\theta\lambda).   
  \end{equation*}
  Then we get
  \begin{align*}
    \Xi(x)\Xi(y)=&\vs_{\lambda+\lambda',\dm}^{-1}q^{\frac{1}{2}(\lambda,\lambda+\theta\lambda)+\frac{1}{2}(\lambda',\lambda'+\theta\lambda')+(\lambda',\lambda+\theta\lambda)}\wObb(xy)\cK_{-\theta(\lambda+\lambda')}\\
    =&\vs_{\lambda+\lambda',\dm}^{-1}q^{\frac{1}{2}(\lambda+\lambda',\lambda+\lambda'+\theta\lambda+\theta\lambda')}\wObb(xy)\cK_{-\theta(\lambda+\lambda')}\\
    =&\Xi(xy).
  \end{align*}
  
For $i\in\wI$, we have $F_i,K_i\tTD_{\bw}^{-1}(E_{\tau i})\in\tU[\alpha_i]$. By \eqref{eq:disting-para} and \eqref{eq:formula of wObb}, we have
\begin{align*}
\Xi(B_i^\sigma)=& \vs_{i,\dm}^{-1}q^{\frac{1}{2}(\alpha_i,\alpha_i+\theta\alpha_i)}\wObb(B_i^\sigma)\cK_{-\theta\alpha_i}\\
=&\vs_{i,\dm}^{-1}q^{\frac{1}{2}(\alpha_i,\alpha_i+\theta\alpha_i)}\wObb(F_i)\cK_{-\theta\alpha_i}+\vs_{i,\dm}^{-1}q^{\frac{1}{2}(\alpha_i,\alpha_i+\theta\alpha_i)}\wObb(K_i\tTD_{\bw}^{-1}(E_{\tau i}))\cK_{-\theta\alpha_i}\\
=&-\vs_{i,\dm}^{-1}q^{\frac{1}{2}(\alpha_i,\alpha_i+\theta\alpha_i)}q_i^{-2}B_i\\
=&B_i.
\end{align*}
The remaining identities in \eqref{eq:formula of Xi on generators}  follow directly from \eqref{eq:formula of wObb}.
\end{proof}

 By \eqref{eq:small weyl group} and \eqref{eq:T retain grade}, if $2\lambda=2\nu-(1+\theta)\mu$, for $i\in\fwItau$, we know
\[2\br_i(\lambda)=2\br_i(\nu)-\br_i((1+\theta)(\mu))=2\br_i(\nu)-(1+\theta)(\br_i(\mu)).\]
Thus, $\tT_{\br_i}^{\pm1}|_{\tU_{\text{int}}}$ is algebra endomorphism of $\tU_{\text{int}}$.
\begin{lemma}\label{lem:Xi commute brT}
  For $i\in\fwItau$, we have $\tT_{\br_i}^{\pm1}\circ\Xi=\Xi\circ\tT_{\br_i}^{\pm1}|_{\tU_{\text{int}}}$.
\end{lemma}
\begin{proof}
    For $x\in\tU[\lambda]\cap\tU_{(\mu,\nu)}$, by \eqref{eq:wObb commute Tbr} and \eqref{eq:formula of wObb}, we have
    \begin{align*}
        \tT_{\br_i}^{\pm1}\circ\Xi(x)=&\vs_{\lambda,\dm}^{-1}q^{\frac{1}{2}(\lambda,\lambda+\theta\lambda)}\cdot\vs^{\frac{1}{2}}_{\lambda+\theta(\br_i\lambda)-\theta\lambda-\br_i\lambda,\dm}\tT_{\br_i}^{\pm1}\circ\wObb(x)\cK_{-\br_i\circ\theta\lambda}\\
        =&\vs_{\lambda,\dm}^{-1}q^{\frac{1}{2}(\lambda,\lambda+\theta\lambda)}\cdot\vs^{\frac{1}{2}}_{\lambda+\theta(\br_i\lambda)-\theta\lambda-\br_i\lambda,\dm}\wObb(\tT_{\br_i}^{\pm1}(x))\cK_{-\theta(\br_i\lambda)}\\
        =&\vs_{\br_i\lambda,\dm}^{-1}q^{\frac{1}{2}(\br_i\lambda,\br_i\lambda+\theta\br_i\lambda)}\wObb(\tT_{\br_i}^{\pm1}(x))\cK_{-\theta(\br_i\lambda)}\\
        =&\Xi\circ\tT_{\br_i}^{\pm1}(x)
    \end{align*}
    The third equality follows from $\vs_{\theta\lambda,\dm}=\vs_{-\lambda,\dm}.$
\end{proof}

\begin{proposition}
  \label{prop:rktwoRPijb}
Let $i\in\fwItau$ and $j\notin \bIi$. Then for the non-commutative polynomials $R_{ij,\bullet}(x_i,x_{\tau i}, y_i, y_{\tau i}, z)$ and $P_{ij,\bullet}(x_i, x_{\tau i}, y_i, y_{\tau i}, z)$ in Corollary \ref{cor:rktwoRPijb}, we have
\begin{align}
  \tT_{\br_i}^{-1} (\tTD_{\bw}( E_{\tau j}) K_j') =& R_{ij,\bullet}(B_i, B_{\tau i}, \cK_i, \cK_{\tau i}, \tTD_{\bw}( E_{\tau j}) K_j'),\label{eq:TF-poly-3}\\
  \tT_{\br_i} (\tTD_{\bw}( E_{\tau j}) K_j') =& P_{ij,\bullet}(\Xi(B_i), \Xi(B_{\tau i}), \cK_i, \cK_{\tau i}, \tTD_{\bw}( E_{\tau j}) K_j')\label{eq:TF-poly-4}.
\end{align}  
\end{proposition}
\begin{proof}
 Applying $\Xi$ to \eqref{eq:TF-poly-1} and \eqref{eq:TF-poly-2},  using Corollary~\ref{cor:rktwoRPijb}, Lemmas \ref{lem:Xi is alge-homomorphism} and \ref{lem:Xi commute brT}, we obtain \eqref{eq:TF-poly-3}, \eqref{eq:TF-poly-4}. 
\end{proof}
Eventually, Corollary~\ref{cor:rktwoRPijb} and Proposition~\ref{prop:rktwoRPijb} together establish \cite[Conjecture 5.13]{WZ23}.

%%%%%%%%%%%%%%%%%%%%%%%
\subsection{New symmetries on $\tUi$}
By the definition of $W^\circ$, it suffices to construct ibraid group symmetries associated to vertices $i\in\fwItau$.
\begin{lemma}\cite[Lemma 5.1]{WZ23}
  \label{lem:rktwo1}
For $i\in\fwItau$ and $j\notin \bIi$, we have
\begin{align}
 F_j  \tfX_i &=\tfX_i  F_j,
 \label{eq:FjUp}\\
\tT_{\br_i}^{-1}\big(\tTD_{\bw}(E_{\tau j})K_j'\big)  \cdot \tfX_i
&= \tfX_i \cdot  \tT_{\br_i}^{-1}\big(\tTD_{\bw}(E_{\tau j}) K_j'\big).
 \label{eq:TEjUp}
\end{align}
\end{lemma}
The following Proposition follows by an argument similar to that of \cite[Theorem 4.7]{WZ23}, using Corollary~\ref{cor:rktwoRPijb}, Proposition~\ref{prop:rktwoRPijb} and Lemma~\ref{lem:rktwo1}.
\begin{proposition}
    \label{prop:newb0}
Let $i\in\fwItau$.
\begin{enumerate}
\item
For any $x^\imath \in \tUi$, there is a unique element $\hat x \in \tUi$ such that $\hat x \tfX_i =\tfX_i \tT_{\br_i}^{-1}(x^\imath)$.
\item
The map $x^\imath \mapsto \hat x$ is an injective endomorphism of the algebra $\tUi$, denoted by $\tTa{i}$.
\end{enumerate}
\end{proposition}

\begin{proposition}[{\cite[Proposition 4.1,Theorem 4.14]{WZ23}\label{prop:iT on bUi}}]
For $i\in\fwItau$, the following formulas hold:
\begin{align*}
    \tTT'_{i,-1}(u)=&\tau_{\bullet,i}\circ\tau(u),\qquad u\in\tbU,\\
    \tTT'_{i,-1}(\tk_j)=&\vs_{\br_i\alpha_j-\alpha_j,\dm}^{-1}\tk_{\br_i\alpha_j},\qquad j\in\wI,\\
    \tTT'_{i,-1}(B_i)=&-q^{(\alpha_i,\theta\alpha_i)}\tTD_{\bw}^2(B_{\tau_{\bullet,i}\tau i})\cK^{-1}_{\tau_{\bullet,i}\tau i}.
\end{align*}
\end{proposition}
Recall the algebra homomorphisms $\Ob$ and $\wObb$ in \eqref{eq:auto-Omega}. Recall the algebra endomorphism $\Xi$ of $\tU_{\text{int}}$ in \eqref{eq:formula of Xi} and the quantum group $\tbiU$ associated to $ \I_{\bullet,i}$. By \eqref{eq:S-sigma-omega=Tbw-tau}, we know 
\begin{equation}\label{eq:local formula of wObb}
    \wObb|_{\tbiU}=\tT_{\br_i}\circ\tTD_{\bw}^2\circ\tau_{\bullet,i}\circ\tau.
\end{equation}
Since $\tT_{\br_i}(\cK^{-1}_{\tau_{\bullet,i}\tau i})=\vs_{i,\dm}^{-2}\cK_{-\theta\alpha_i}$, using Proposition~\ref{prop:iT on bUi}, we have
\begin{equation*}
 \tT_{\br_i}(\tTT'_{i,-1}(B_i))=-q_i^2 \wObb(B_i)\cK_{-\theta\alpha_i}=\Xi(B_i).   
\end{equation*}
So if we write $\widehat{B}_i:=\tT_{\br_i}(\tTT'_{i,-1}(B_i))$, by Proposition~\ref{prop:rktwoRPijb}, \eqref{eq:TF-poly-4} can be reformulated as 
\begin{equation}\label{eq:TF-poly-5}
 \tT_{\br_i} (\tTD_{\bw}( E_{\tau j}) K_j') = P_{ij,\bullet}(\widehat{B}_i, \widehat{B}_{\tau i}, \cK_i, \cK_{\tau i}, \tTD_{\bw}( E_{\tau j}) K_j')   
\end{equation}

    \begin{theorem}
		\label{thm: relative T_i conjugate}
		For $i\in\fwItau$, there are mutually inverse automorphisms $\tTT'_{i,-1}$ and $\tTT_{i,+1}''$ on $\tUi$ such that
		\begin{align}
			\label{eq:relbraid1}
			\tTT'_{i,-1}(x)\widetilde{\Upsilon}_i &=\widetilde{\Upsilon}_i\tT_{\br_i}^{-1}(x),
			\\
			\label{eq:relbraid2}
			\tTT_{i,+1}''(x)\tT_{\br_i}(\widetilde{\Upsilon}_i)^{-1} &=\tT_{\br_i}(\widetilde{\Upsilon}_i)^{-1}\tT_{\br_i}(x).
		\end{align}
		Moreover, $\tTT'_{i,-1}=\sigma^\imath\circ \tTT_{i,+1}''\circ \sigma^\imath$, and there exists a group homomorphism $\mathrm{Br}(W^\circ)\rightarrow \Aut(\tUi)$, $\br_i\mapsto \tTT''_{i,+1}$ for $i\in\fwItau$.
	\end{theorem}
\begin{proof}
Equation \eqref{eq:relbraid1} follows from Proposition~\ref{prop:newb0}.  \eqref{eq:relbraid2} follows by adapting the  proof of \cite[Theorem 6.1]{WZ23}, using \eqref{eq:TF-poly-2} and \eqref{eq:TF-poly-5}. The remaining statements are followed by the same arguments as in \cite[\S6.5, \S9.1]{WZ23} and \cite[\S3.4]{Z23}.
\end{proof}

%\red{No need to give short hand for $\tTT'_{i,-1}$, $\tTT''_{i,+1}$?}

%%%%%%%%%%%%%%%%%%%%%
\section{The tensor quasi \texorpdfstring{$K$}{K}-matrix and intertwining properties}
\label{sec:tensor K}

This section uses the tensor quasi $K$-matrix to relate the braid group symmetries on quantum groups and iquantum groups.

\subsection{The tensor quasi \texorpdfstring{$K$}{K}-matrix}
Set $\tU^{\ge0}=\tU^0 \tU^+$ and $\tU^{\le0}=\tU^-\tU^0$. These are Hopf subalgebras of $\tU$. 
\begin{proposition}[{\cite[Proposition 1.2.3]{Lus93}}]\label{bilinear form}
     There exists a bilinear pairing $\langle \ ,\ \rangle : \tU^{\leq 0} \times \tU^{\geq 0} \rightarrow \K$ such that for all $x, x' \in \tU^{\geq 0}$, all $y, y' \in \tU^{\leq 0}$, all $\mu, \nu \in \N\I$, and all $i,j\in \I$
     \begin{align}
         \langle y, xx'\rangle = \langle\Delta(y), x' \otimes x\rangle,& \quad \langle yy', x\rangle = \langle y \otimes y', \Delta(x)\rangle,\label{eq:form-kow}\\
         \langle K_\mu, K_\nu\rangle = q^{-(\mu, \nu)}, \quad \langle K_\mu', K_\nu\rangle =&\langle K_\mu, K_\nu'\rangle= q^{(\mu, \nu)},\quad \langle K_\mu', K_\nu'\rangle = q^{-(\mu, \nu)},\\
         \langle F_i, E_j\rangle =& \delta_{ij}\frac{1}{q_i^{-1}-q_i},\\
         \langle K_\mu, E_j\rangle=\langle K_\mu', E_j\rangle = 0, &\quad \langle F_i, K_\mu\rangle =\langle F_i, K_\mu'\rangle= 0.
     \end{align}
\end{proposition}

This pairing $\langle \ ,\ \rangle : \tU^{\leq 0} \times \tU^{\geq 0} \rightarrow \K$ is called a skew-Hopf pairing of $\tU$. Note that the pairing $\langle \ ,\ \rangle$ restricted to $\tU^- \times \tU^+$ is nondegenerate.

For $\mu\in \N\I$, set $r(\mu)=\dim(\tU^+_\mu)$ and let $\{u^+_{\mu,i}\,|\,i=1,\dots,r(\mu)\}$ be a basis of $\tU^+_\mu$ with dual basis $\{u^-_{\mu,i}\,|\,i=1,\dots,r(\mu)\}\subset \tU^-_{-\mu}$. Define $\Theta_\mu=\sum_{i=1}^{r(\mu)} u^-_{\mu,i}\otimes u^+_{\mu,i}$. We denote the completions of the algebras $\tU\otimes\tU$ and $\tU$ by $\tU \widehat{\otimes} \tU$ and $\widehat{\tU}$, respectively. The quasi $R$-matrix of $\tU$ is defined by
\begin{align}\label{eq:Theta-def}
  \Theta=\sum_{\mu\in \N\I} \Theta_\mu\in\tU \widehat{\otimes} \tU,
\end{align}
which is independent of the bases.
By \cite[Theorem 4.1.2]{Lus93}, the quasi $R$-matrix satisfies the relation
\begin{align}\label{eq:Theta-intertwiner}
  \Delta(\overline{u})\cdot \Theta=\Theta\cdot (\obar \otimes \obar)\circ\Delta(u) \qquad \text{for all $u\in \tU$.}
\end{align}  
 Define $\overline{\Theta}=(\obar\otimes\obar)(\Theta)$. By construction, $\Theta$ is invertible in $\tU \widehat{\otimes} \tU$, and $\Theta$ is uniquely determined by the intertwiner property \eqref{eq:Theta-intertwiner} and the normalization $\Theta_0=1\otimes 1$. These two properties together imply that $\Theta^{-1}=\overline{\Theta}$.

We now return to the Satake diagram $(\I=\I_\bullet\cup\I_\circ,\tau)$. Let $i\in\fwItau$. For the subset $\bI\subset\I$ (resp. $\bIi\subset\I$), let $\Theta_\bullet$ (resp. $\Theta_{\bullet,i}$)  denote the quasi $R$-matrix of $\tbU$ (resp. $\tbiU$) and denote its inverse by $\ov{\Theta}_\bullet=(\obar\otimes\obar)(\Theta_\bullet)$ (resp. $\ov\Theta_{\bullet,i}$). Write $\cR_\bullet^{+,l}:=\omega(\cR_\bullet^{-,r})$.

For any $\mu\in \N\I$, set $s(\mu)=\dim((\cR_\bullet^{+,l})_\mu)=\dim((\cR_\bullet^{-,r})_{-\mu})$ and let $\{e_{\mu,i}\,|\,i=1\dots,s(\mu)\}$ be a basis of $(\cR_\bullet^{+,l})_\mu$ with dual basis $\{f_{\mu,i}\,|\,i=1\dots,s(\mu)\}\subset(\cR_\bullet^{-,r})_{-\mu}$ for the pairing $\langle\ ,\ \rangle$. Define $\Lambda_{\bullet,\mu}=\sum_{i=1}^{s(\mu)} f_{\mu,i}\otimes e_{\mu,i}$ and set
\begin{align*}
  \Lambda_\bullet=\sum_{\mu\in \N\I} \Lambda_{\bullet,\mu}\in \tU \widehat{\otimes} \tU;
\end{align*} 
see \cite[\S 5.2]{KY26}.

\begin{lemma}[{\cite[Lemma 5.4]{KY26}}]\label{lem:rel-quasi} The following statements hold.
  \begin{enumerate}
    \item $\Lambda_\bullet=\Theta \overline{\Theta}_\bullet$.
    \item For all $h\in \tbU\tU^0$, the relation $\Delta(h)\Lambda_\bullet=\Lambda_\bullet \Delta(h)$
      holds in $\tU \widehat{\otimes} \tU$.
  \end{enumerate}  
\end{lemma}

Recall that the Letzter map $\psi:\tUi\rightarrow \cA$ is a linear isomorphism, and set
\begin{align}\label{eq:TB-def}
  \Theta^\imath=(\psi^{-1}\otimes \id)(\Theta) \in \prod_{\mu\in \N\I} \tUi\otimes \tU^+_\mu.
\end{align}
We call $\Theta^\imath$ the tensor quasi $K$-matrix. 

The map $\sigma_\tau=\sigma^\imath\circ\tau:\tUi\rightarrow \tUi$ is an anti-automorphism. 
\begin{proposition}[{\cite[Theorem 5.8]{KY26}}]
    \label{thm:intertwiner2}
  For $b\in \tUi$, the following identity holds in $\tU \widehat{\otimes} \tU$:
  \begin{align}\label{eq:intertwiner2}
    \Delta(\sigma_\tau(b)) \cdot\Theta^\imath = \Theta^\imath\cdot (\sigma_\tau\otimes(\sigma\circ\tau))\circ \Delta(b).
  \end{align}
\end{proposition}

\begin{proposition}[{\cite[Corollaries 5.10, 5.11]{KY26}}]
The following identities hold:
 \begin{equation}\label{description of the tqK}
     (\varepsilon\otimes 1)(\Theta^\imath)=\tfX,\qquad \Theta^\imath= \Delta(\tfX)\cdot \Theta\cdot(\tfX^{-1}\otimes 1).
 \end{equation}
\end{proposition}
    
\subsection{The tensor intertwining relation}
Recall that $\Theta_{\bullet,i}$ and $\Theta_\bullet$ are the ordinary quasi
$R$-matrices of $\tU_{\bullet,i}$ and $\tbU$, respectively. Set
\begin{align}
    \Theta_{\bullet,i}^{\imath}
&=(\psi^{-1}\otimes\id)(\Theta_{\bullet,i}),
\qquad\Lambda_{\bullet,i}=\Theta_{\bullet,i}\overline\Theta_\bullet.
\label{eq:local tensor-definition}
\end{align}

Remark~\ref{rem:local Letzter and Letzter} and equation~\eqref{description of the tqK} give
\begin{equation}\label{local-tensor to K}
    (\varepsilon\otimes 1)(\Theta_{\bullet,i}^{\imath})=\tfX_i.
\end{equation}
\begin{lemma}
\label{lem:tensor-factorization}
The following factorization holds:
\begin{align}
    \Theta_{\bullet,i}^{\imath}\overline\Theta_\bullet=\Delta(\tfX_i)\Lambda_{\bullet,i}(\tfX_i^{-1}\otimes1)=(\psi^{-1}\otimes\id)(\Lambda_{\bullet,i}).\label{eq:tensor-factorization}
\end{align}
\end{lemma}
\begin{proof}
Equations~\eqref{eq:fX2} and \eqref{description of the tqK} give 
\begin{align*}
    \Theta_{\bullet,i}^{\imath}\overline\Theta_\bullet&=\Delta(\tfX_i)\cdot \Theta_{\bullet,i}\cdot(\tfX_i^{-1}\otimes 1)\ov{\Theta}_\bullet\\
    &=\Delta(\tfX_i)\cdot \Theta_{\bullet,i}\ov{\Theta}_\bullet(\tfX_i^{-1}\otimes 1)
    \\
    &=\Delta(\tfX_i)\Lambda_{\bullet,i}(\tfX_i^{-1}\otimes1).
\end{align*}
The $\cH_\bullet$-bimodule property of $\psi$ gives 
\begin{align*}
    (\psi^{-1}\otimes\id)(\Lambda_{\bullet,i})=&(\psi^{-1}\otimes\id)(\Theta_{\bullet,i}\ov{\Theta}_\bullet)
    \\
    =&(\psi^{-1}\otimes\id)(\Theta_{\bullet,i})\ov{\Theta}_\bullet
    \\
    =&\Theta_{\bullet,i}^{\imath}\overline\Theta_\bullet.
\end{align*}
\end{proof}

\begin{proposition}
\label{lem:mixed-coproduct}
For any $u\in\tU$, one has
\begin{align}
\Delta(\tTD'_{\br_i,-1}(u))\Lambda_{\bullet,i}
&=\Lambda_{\bullet,i}
(\tTD'_{\br_i,-1}\otimes\tTD'_{\br_i,-1})\Delta(u),
\label{eq:ordinary-coproduct}\\
\Delta(\tT'_{\br_i,-1}(u))\Lambda_{\bullet,i}
&=\Lambda_{\bullet,i}
(\tT'_{\br_i,-1}\otimes\tTD'_{\br_i,-1})\Delta(u).
\label{eq:mixed-coproduct}
\end{align}
\end{proposition}
\begin{proof}
    We first recall the ordinary braid coproduct formula
\cite[37.3.2]{Lus93}, in the conventions of \cite[(1.2), \S4.4]{WZ23}:
\begin{align}
\Delta(\tTD'_{j,-1}(u))\mathbf L_j''
=\mathbf L_j''(\tTD'_{j,-1}\otimes\tTD'_{j,-1})\Delta(u),
\label{eq:simple-coproduct}
\end{align}
where $\mathbf L_j''=\Theta_j$ is the quasi
$R$-matrix, explicitly given by
\begin{align}
\mathbf L_j''=\sum_{n\geq0}(-1)^nq_j^{-n(n-1)/2}
(q_j-q_j^{-1})^n[n]_j^!\,F_j^{(n)}\otimes E_j^{(n)}.
\label{eq:elementary-R}
\end{align}
For $w\in W_{\bullet,i}$, choose a reduced expression $w=s_{i_1}\cdots s_{i_m}$ and set $\mathbf{i}=(i_1,\ldots,i_m)$. Define $\mathbf L_{\bi}''$ by
\begin{equation}\label{def:R-Product}
  \mathbf L_\bi''=\prod_{k=1}^m \Big(\tTD'_{s_{i_1}\cdots s_{i_{k-1}},-1}\otimes\tTD'_{s_{i_1}\cdots s_{i_{k-1}},-1}\Big)(\mathbf L_{i_k}'')  
\end{equation}
and set $\mathbf L_\varnothing''=1\otimes1$. The PBW theorem gives
\begin{equation*}
    \mathbf L_\bi''\in\prod_{\mu\in\N\I}\tU^-_{-\mu}\otimes\tU^+_{\mu},\qquad (\mathbf L_w'')_0=1\otimes1.
\end{equation*}
Equation~\eqref{eq:simple-coproduct} yields
\begin{align}
    \Delta(\tTD'_{w,-1}(u))\mathbf L_\bi''
&=\mathbf L_\bi''(\tTD'_{w,-1}\otimes\tTD'_{w,-1})\Delta(u).
\label{eq:pre-iterated-coproduct}
\end{align}
Lemma~\ref{lem:Tw-on-X} and \eqref{eq:pre-iterated-coproduct} imply $\mathbf L_{\bi}''=\Theta_{\bullet,i}$ (resp. $\mathbf L_{\bj}''=\Theta_\bullet$) for any reduced expression $\bi$ (resp. $\bj$) of $\bwi$ (resp. $\bw$).

Assume $\bi$ and $\bj$ are reduced expressions of $w$ and write $v=w^{-1}\bwi$. Fix a reduced expression $\bk$ of $v$. Then $\bi\bk$ and $\bj\bk$ are reduced expressions of $\bwi$. Hence
\begin{equation*}
  \Theta_{\bullet,i}=\mathbf L_\bi''(\tTD'_{w,-1}\otimes\tTD'_{w,-1})(\mathbf L_{\bk}'')=\mathbf L_\bj''(\tTD'_{w,-1}\otimes\tTD'_{w,-1})(\mathbf L_{\bk}'').  
\end{equation*}
It follows that \eqref{def:R-Product} is independent of the choice of reduced expression. Thus, we may write $\mathbf L_{w}'':=\mathbf L_{\bi}''$ for any reduced expression $\bi$ of $w$. Equations~\eqref{eq:elementary-R} and \eqref{def:R-Product} give
\begin{align}
\Delta(\tTD'_{w,-1}(u))\mathbf L_w''
&=\mathbf L_w''(\tTD'_{w,-1}\otimes\tTD'_{w,-1})\Delta(u),
\label{eq:iterated-coproduct}\\
\mathbf L_{vw}''
&=\mathbf L_v''(\tTD'_{v,-1}\otimes\tTD'_{v,-1})(\mathbf L_w'')
\quad\text{if }\ell(vw)=\ell(v)+\ell(w).
\label{eq:R-product}
\end{align}
Moreover, $\tTD'_{\br_i,-1}$ restricts to the diagram automorphism
$\tau_{\bullet,i}\circ\tau$ of $\tbU$. Applying
\eqref{eq:R-product} to \eqref{def:bsi} gives
\begin{equation*}
    \Theta_{\bullet,i}
=\mathbf L_{\br_i}''
(\tTD'_{\br_i,-1}\otimes\tTD'_{\br_i,-1})(\Theta_\bullet)
=\mathbf L_{\br_i}''\Theta_\bullet.
\end{equation*}
Lemma~\ref{lem:rel-quasi} therefore gives $\mathbf L_{\br_i}''=\Lambda_{\bullet,i}$, while
\eqref{eq:iterated-coproduct} proves \eqref{eq:ordinary-coproduct}.

Observe that
\begin{align}
(\tPsi_{\bvs_\dm}^{-1}\otimes\id)(\Lambda_{\bullet,i})
=\Lambda_{\bullet,i}.
\label{eq:scaling-R}
\end{align}
A direct check on generators gives
\begin{align}
\Delta\circ\tPsi_{\bvs_\dm}^{\pm1}
=(\tPsi_{\bvs_\dm}^{\pm1}\otimes\id)\circ\Delta.
\label{eq:scaling-coproduct}
\end{align}
This proves \eqref{eq:mixed-coproduct}. 
\end{proof}

\iffalse
\begin{remark}
The second operator in \eqref{eq:mixed-coproduct} is
$\tTD'_{\br_i,-1}$, not $\tT'_{\br_i,-1}$: the scaling identity
\eqref{eq:scaling-coproduct} acts on the first tensor factor only.
Taking $\varepsilon\otimes\id$ in \eqref{eq:mixed-coproduct} also yields
\begin{align}
(\varepsilon\circ\tT'_{\br_i,-1}\otimes\tTD'_{\br_i,-1})\Delta(u)
=\tT'_{\br_i,-1}(u),
\label{eq:mixed-counit}
\end{align}
since $(\varepsilon\otimes\id)(\Lambda_{\bullet,i})=1$.
\end{remark}
\fi 

\begin{theorem}
\label{thm:tensor}
 In $\tU\widehat{\otimes}\tU$, the following identity holds:
\begin{align}
\Delta\bigl(\tTa{i}(x)\bigr)
\Theta_{\bullet,i}^{\imath}\overline\Theta_\bullet=&\Theta_{\bullet,i}^{\imath}\overline\Theta_\bullet
(\tTa{i}\otimes\tTD'_{\br_i,-1})\Delta(x),
\qquad x\in\tUi.
\label{eq:tensor-braid}
\end{align}
\end{theorem}

\begin{proof}
Since $\tfX_i$ is invertible, \eqref{eq:relbraid1} gives
\begin{align}
\tTa{i}(x)=\tfX_i\tT'_{\br_i,-1}(x)\tfX_i^{-1},
\qquad x\in\tUi.
\label{eq:conjugation}
\end{align}
Combining Lemma \ref{lem:tensor-factorization} and Proposition \ref{lem:mixed-coproduct} gives
\begin{align*}
\Delta\bigl(\tTa{i}(x)\bigr)
\Theta_{\bullet,i}^{\imath}\overline\Theta_\bullet
&\quad=\Delta(\tfX_i)\Delta(\tT'_{\br_i,-1}(x))
\Lambda_{\bullet,i}(\tfX_i^{-1}\otimes1)\\
&\quad=\Delta(\tfX_i)\Lambda_{\bullet,i}
(\tT'_{\br_i,-1}\otimes\tTD'_{\br_i,-1})\Delta(x)
(\tfX_i^{-1}\otimes1)\\
&\quad=\Theta_{\bullet,i}^{\imath}\overline\Theta_\bullet
(\tfX_i\otimes1)
(\tT'_{\br_i,-1}\otimes\tTD'_{\br_i,-1})\Delta(x)
(\tfX_i^{-1}\otimes1)\\
&\quad=\Theta_{\bullet,i}^{\imath}\overline\Theta_\bullet
(\tTa{i}\otimes\tTD'_{\br_i,-1})\Delta(x).
\end{align*}
The last equality holds because
$\Delta(\tUi)\subseteq\tUi\otimes\tU$, so
\eqref{eq:conjugation} applies to the first tensor components of
$\Delta(x)$. This proves \eqref{eq:tensor-braid}.
\iffalse
Conversely, assume \eqref{eq:tensor-braid} and
\eqref{eq:counit-normalization}. Lemma
\ref{lem:tensor-factorization} implies
\begin{align}
(\varepsilon\otimes\id)
(\Theta_{\bullet,i}^{\imath}\overline\Theta_\bullet)=\tfX_i.
\label{eq:relative-counit}
\end{align}
Applying $\varepsilon\otimes\id$ to \eqref{eq:tensor-braid} and then using \eqref{eq:counit-normalization} and \eqref{eq:mixed-counit} gives
\begin{align*}
\tTa{i}(x)\tfX_i
&=\tfX_i(\varepsilon\circ\tTa{i}\otimes
\tTD'_{\br_i,-1})\Delta(x)\\
&=\tfX_i(\varepsilon\circ\tT'_{\br_i,-1}\otimes
\tTD'_{\br_i,-1})\Delta(x)\\
&=\tfX_i\tT'_{\br_i,-1}(x).
\end{align*}
\fi
\end{proof}

%%%%%%%%%%%%%%%%%%
\section{\texorpdfstring{$\mathrm i$}{i}Braid group symmetries on \texorpdfstring{$\cR_{\bullet,i}^{-,r}$}{R}}
\label{sec:R}

Throughout this section, $\psi:\tU^\poly\rightarrow \cA$ denotes the natural projection; see \eqref{eq:Upoly-decomp}. Let $i\in\fwItau$ and write $\J:=\I-\bIi$. %The assumption will be used consistently in all subsequent arguments.
\subsection{The Letzter map on parabolic subalgebras}
Following Subsection~\ref{sec:letzter map}, define
\begin{align}\label{eq:generator of radicals}
\cR_{\bullet,i}^{-,r}&=\langle\ad_r(\tbiU^-)(F_j)\mid j\in\J\rangle,\qquad \cR_i= \langle\ad_r(\tbU^-)(F_i),\ad_r(\tbU^-)(F_{\tau i})\rangle\end{align}
   as subalgebras of $\tU$.
   Thus,
   \begin{align}
  \cR_{\bullet,i}^{-,r}  &=\tU^-\cap\tTD_{\bwi}(\tU^-)=\bigcap_{j\in\bIi}(\ker\partial_j^R|_{\tU^-}),\label{eq:radicals} \\
    \cR_i&=\tbiU^-\cap\tTD_{\bw}(\tbiU^-)=\bigcap_{j\in\bI}(\ker\partial_j^R|_{\tbiU^-})
\label{eq:local-radicals}.
\end{align}
\iffalse
\begin{align}
    \cR_{\bullet,i}^{-,r}&=\tU^-\cap\tTD_{\bwi}(\tU^-)=\langle\ad_r(\tbiU^-)(F_j)\mid j\in\J\rangle=\bigcap_{j\in\J}(\ker\partial_j^R|_{\tU^-})\label{eq:radicals},\\
   \begin{split} \cR_i&=\tbiU^-\cap\tTD_{\bw}(\tbiU^-)=\langle\cup_{j\in\bI}\ad_r(\tbU^-)(F_i),\cup_{j\in\bI}\ad_r(\tbU^-)(F_{\tau i})\rangle=\bigcap_{j\in\{i,\tau i\}}(\ker\partial_j^R|_{\tbiU^-})
   \end{split}\label{eq:local-radicals}.
\end{align}
\fi
The algebras $\cR_{\bullet,i}^{-,r}$ and $\cR_i$ are right coideal subalgebras. As in Subsection \ref{sec:letzter map}, multiplication gives
\begin{align}
    \tbiU^-\cong\cR_i\otimes\tbU^-,\quad\cR_\bullet^{-,r}\cong\cR_{\bullet,i}^{-,r}\otimes\cR_i,\quad \cA\cong\cR_{\bullet,i}^{-,r}\otimes\cR_i\otimes\cH_\bullet.
\label{eq:decomposition}
\end{align}

\begin{lemma}\label{lem:support}
For $u\in\tU^-$ and $x\in\cR_{\bullet,i}^{-,r}$, respectively,
\begin{align}
\psi^{-1}(u)-u
&\in\sum_{0\ne\beta\in\N\I}\tU^-\tU^+K_\beta',
\label{eq:support-all}\\
\psi^{-1}(x)-x
&\in\sum_{\beta\in\N\I\setminus\N \bIi}\tU^-\tU^+K_\beta'.
\label{eq:support-parabolic}
\end{align}
\end{lemma}
\begin{proof}
 Set $\bDelta=(\overline{\phantom{x}}\otimes\overline{\phantom{x}})
\circ\Delta\circ\overline{\phantom{x}}$. Equations~\eqref{description of the tqK} and \eqref{eq:Theta-intertwiner} give
\begin{align}
    \Theta^\imath=(\psi^{-1}\otimes \id)(\Theta)=\Theta\bDelta(\tfX)(\tfX^{-1}\otimes1).
\label{eq:tensor-expand}
\end{align}
Set $C:=\bDelta(\tfX)(\tfX^{-1}\otimes1)$. Since
\begin{align*}
\bDelta(\tU^+_\eta)
&\subseteq\sum_{\beta+\gamma=\eta}
\tU^+_\gamma K_\beta'\otimes\tU^+_\beta,
& (\tU^+K_\beta')\tU^+&\subseteq\tU^+K_\beta',
\end{align*}
we may decompose $C=\sum_{\beta\in\N\I}C_\beta$ by the weight of the
second tensor factor. More precisely,
\begin{align}
 C_0=1\otimes 1,\qquad C_\beta=\sum_{s_\beta} c_{\beta,s_\beta}K_\beta'\otimes v^+_{\beta,s_\beta}\in\widehat{\tU^+}K_\beta'
\widehat\otimes\tU^+_\beta.
\label{eq:dual-bases}   
\end{align}
Choose dual bases for $\tU^\pm$ so that
\begin{align*}
   \Theta=\sum_{\mu,r_\mu}u^-_{\mu,r_\mu}\otimes u^+_{\mu,r_\mu},
\qquad\langle u^-_{\mu,r_\mu},u^+_{\nu,r_\nu}\rangle=\delta_{\mu,\nu}\delta_{r_\mu ,r_\nu}. 
\end{align*}
For $u\in\tU^-$, define $\ev_u(v)=\langle u,v\rangle$ for $v\in\tU^+$. Equation~\eqref{eq:tensor-expand} then gives
\begin{align*}
(\id\otimes\ev_u)
\bigl((\psi^{-1}\otimes\id)(\Theta)\bigr)
&=\sum_{\mu,r_\mu}\psi^{-1}(u^-_{\mu,r_\mu})
\langle u,u^+_{\mu,r_\mu}\rangle\nonumber\\
&=\psi^{-1}\Bigl(\sum_{\mu,r_\mu}
u^-_{\mu,r_\mu}\langle u,u^+_{\mu,r_\mu}\rangle\Bigr)
=\psi^{-1}(u),%\label{eq:pair-left}
\\
(\id\otimes\ev_u)(\Theta C)
&=u+\sum_{\substack{\beta\ne0,\,\mu\\r_\mu,s_\beta}}
u^-_{\mu,r_\mu}c_{\beta,s_\beta}
\langle u,u^+_{\mu,r_\mu}v^+_{\beta,s_\beta}\rangle K_\beta'.
%\label{eq:pair-right}
\end{align*}
Consequently,
\begin{align}
\psi^{-1}(u)=u+\sum_{\substack{\beta\ne0,\,\mu\\r_\mu,s_\beta}}
u^-_{\mu,r_\mu}c_{\beta,s_\beta}
\langle u,u^+_{\mu,r_\mu}v^+_{\beta,s_\beta}\rangle K_\beta',
\label{eq:inverse-expansion}
\end{align}
which proves \eqref{eq:support-all}.

Now let $x\in\cR_{\bullet,i}^{-,r}$ and write
$\Delta(x)=\sum_m x_m\otimes y_m$ with homogeneous
$x_m\in\cR_{\bullet,i}^{-,r}$. The skew-Hopf pairing gives
\begin{align}
\langle x,u^+_{\mu,r_\mu}v^+_{\beta,s_\beta}\rangle
&=\sum_m\langle x_m,v^+_{\beta,s_\beta}\rangle
\langle y_m,u^+_{\mu,r_\mu}\rangle.\label{eq:pair-order}
\end{align}
For $0\ne\beta\in\N \bIi$, \eqref{eq:generator of radicals} implies $\langle x_m,v^+_{\beta,s_\beta}\rangle=0$.
Thus only $\beta\notin\N \bIi$ remain in
\eqref{eq:inverse-expansion}. 
\end{proof}
%%%%%%%%%%%%%%%%%
\subsection{Specialization and Cartan support of the star product}
Define
\begin{equation}\label{eq:local}
    \cH_{\bullet,i}=\langle\tbU,\tk_i^{\pm1},\tk_{\tau i}^{\pm1}\rangle,\qquad
\cA^{[0]}=\tU^-\cH_{\bullet,i}
=\cR_{\bullet,i}^{-,r}\cR_i\cH_{\bullet,i}.
\end{equation}
Thus,
\begin{align}
\cA&=\bigoplus_{\nu\in\Z\J}
\cA^{[0]}\tk_{\nu},
&\pi_{\bullet,i}\Bigl(\sum_{\nu}
a_{\nu}\tk_{\nu}\Bigr)&=a_{0}, \quad a_\nu\in\cA^{[0]}.
\label{eq:pJ}
\end{align}
Here $\pi_{\bullet,i}:\cA\to\cA^{[0]}$ is a linear projection with respect to the direct sum decomposition \eqref{eq:pJ}. Let $d_{\bullet,i}$ be the function defined on $\Z\I$ by
\begin{align}
d_{\bullet,i}\Bigl(\sum_jm_j\alpha_j\Bigr)&=\sum_{j\in\J}m_j.
\end{align}
Accordingly, set 
\begin{align}
&\cR_{\bullet,i}^{-,r}[n]
=\bigoplus_{d_{\bullet,i}(\mu)=n}(\cR_{\bullet,i}^{-,r})_{-\mu},\quad \text{ for any }n\in\N.
\label{eq:outside-degree}
\end{align}
For convenience, also set $\cR_{\bullet,i}^{-,r}[n]=0$ for $n<0$.

Let $d_j=d_{\bullet,i}(\alpha_j)$, so that $d_j=0$ for $j\in\bIi$ and $d_j=1$ for $j\in\J$. Define a $\Z$-grading on $\tU$ by
\begin{align}
    \deg_{\bullet,i}^+ F_j=0,\quad\deg_{\bullet,i}^+ E_j=d_j,\quad\deg_{\bullet,i}^+ K_j^{\pm1}=\deg_{\bullet,i}^+ K_j'^{\pm1}=\pm d_j.\label{eq:crucial-deg}
\end{align}
For an indeterminate $t$ over $\K$, define an automorphism $\rho_t$ of $\tU\otimes_\K\K[t,t^{-1}]$ by
\begin{align}
\begin{split}
\rho_t(F_j)&=F_j,\;\;\qquad\qquad\quad
\rho_t(E_j)=t^{d_j}E_j,\\
\rho_t(K_j^{\pm1})&=t^{\pm d_j}K_j^{\pm1},\qquad
\rho_t(K_j'^{\pm1})=t^{\pm d_j}K_j'^{\pm1}.
\label{eq:rho}
\end{split}
\end{align}
A direct check on generators shows that $\tU^\poly$ and $\cA$ are graded subalgebras and that $\cJ_\bullet$ is a graded left ideal of $\tU^\poly$. Since $\tU^\poly,\cA,\cJ_\bullet$ are graded subspaces of $\tU$, the automorphism $\rho_t$ preserves $\tU^\poly,\cA,\cJ_\bullet$ over $\K[t,t^{-1}]$.

$\tUi$ itself is not graded. Since $\psi:\tU^\poly\rightarrow \cA$ is the projection onto a graded summand, it is a graded map. Over $\K[t,t^{-1}]$, this gives
\begin{align}
\psi\circ\rho_t=\rho_t\circ\psi,\qquad
\rho_t(h)=h\quad\text{for $h\in\cA^{[0]}$}.
\label{eq:rho-psi}
\end{align}

\begin{lemma}
  For $h\in\cA^{[0]}$ and $z\in\cR_{\bullet,i}^{-,r}$, the following hold: 
\begin{align}
\rho_t(\psi^{-1}(h))\in\tU^\poly[t],\qquad
\rho_t(\psi^{-1}(z))|_{t=0}=z.
\label{eq:polynomial and constant-z}
\end{align}
\end{lemma}
\begin{proof}
 By \eqref{eq:local}, write $h=\sum_lu_lh_l$ for $u_l\in\tU^-$ and $h_l\in\cH_{\bullet,i}$. Since $\psi$ is an $\cH_\bullet$-bimodule homomorphism, \eqref{eq:support-all} and \eqref{eq:rho} give
 \begin{equation*}
  \rho_t(\psi^{-1}(h))=\sum_l\rho_t(\psi^{-1}(u_l))h_l\in\tU[t]\cap\tU^\poly[t,t^{-1}]\subset\tU^\poly[t].   
 \end{equation*}
 The preceding argument also gives $\rho_t(\psi^{-1}(z))\in\tU^\poly[t]$. Equation~\eqref{eq:support-parabolic} further gives
 \begin{equation*}
   \rho_t(\psi^{-1}(z))-z\in t\tU[t].  
 \end{equation*}
 Hence $\rho_t(\psi^{-1}(z))|_{t=0}=z$.
\end{proof}

\begin{proposition}\label{prop:star}
For $x\in\cR_{\bullet,i}^{-,r}[m]$ and
$y\in\cR_{\bullet,i}^{-,r}[n]$, one has
\begin{align}
x*y&=xy+\sum_{0\neq\nu\in\N\J}
z_{\nu}\tk_\nu,\quad z_{\nu}\in\cR_{\bullet,i}^{-,r}[m+n-2d_{\bullet,i}(\nu)]\cR_i\cH_{\bullet,i}.
\label{eq:triangularity}
\end{align}
In particular,
\begin{align}
\pi_{\bullet,i}(x*y)=xy,\quad
x*y-xy\in\bigoplus_{k=1}^{\lfloor(m+n)/2\rfloor}
\cR_{\bullet,i}^{-,r}[m+n-2k]\cR_i\cH_\bullet.
\label{eq:star-projection}
\end{align}
\end{proposition}
\begin{proof}
The defining relations of $\tU$, together with Lemma~\ref{lem:support}, give
\begin{align}
\psi^{-1}(u)\psi^{-1}(v)
\in\sum_{\lambda,\mu\in\N\I}
\tU^-\tU^+K_\lambda K_\mu'
\qquad u,v\in\tU^-.
\label{eq:nonnegative-PBW}
\end{align}
Equations~\eqref{eq:J-triang} and \eqref{eq:nonnegative-PBW} imply
\begin{align}
u*v=\psi(\psi^{-1}(u)\psi^{-1}(v))\in\sum_{\nu\in\N\J}
\cA^{[0]}\tk_\nu.
\label{eq:nonnegative-star}
\end{align}

Equations~\eqref{eq:rho-psi} and \eqref{eq:polynomial and constant-z} yield
\begin{align}
\rho_t(x*y)
=\psi\bigl(\rho_t(\psi^{-1}(x))
\rho_t(\psi^{-1}(y))\bigr),\qquad
\left.\rho_t(x*y)\right|_{t=0}=\psi(xy)=xy.
\label{eq:rho-star}
\end{align}
It follows from \eqref{eq:nonnegative-star} that
\begin{align}
\pi_{\bullet,i}(x*y)=\left.\rho_t(x*y)\right|_{t=0}=xy.
\label{eq:constant-term}
\end{align}

To determine the degrees of the remaining terms, we use another $\Z$-grading of $\tU$ given by
\begin{align}
\deg_{\bullet,i}^- F_j&=d_j,
&\deg_{\bullet,i}^- E_j&=0,
&\deg_{\bullet,i}^- K_j^{\pm1}&=\deg_{\bullet,i}^- K_j'^{\pm1}=\pm d_j.
\label{eq:auxiliary}
\end{align}
Hence $\tU^\poly,\cA$, and $\tUi$ are graded subalgebras.
The kernel $\cJ_\bullet$ is graded as well.
For a graded subspace $V\subseteq\tU$, write
$V\langle r\rangle=V\cap\tU\langle r\rangle$, where
$\tU\langle r\rangle$ denotes the subspace of degree $r$. Hence
\begin{align}
\tU^\poly\langle r\rangle
&=\cA\langle r\rangle\oplus\cJ_\bullet\langle r\rangle,
&\psi(\tUi\langle r\rangle)&\subseteq\cA\langle r\rangle.
\label{eq:graded-projection}
\end{align}
To check that $\psi^{-1}:\cA\rightarrow\tUi$ is a graded map, let $a\in\cA\langle r\rangle$ and write
$\psi^{-1}(a)=\sum_s b_s$ where $b_s\in\tUi\langle s\rangle$. Hence $a=\sum_s\psi(b_s)$ and $\psi(b_s)\in\cA\langle s\rangle$. Thus $\psi(b_s)=0$ for $s\ne r$. Since $\psi:\tUi\rightarrow\cA$ is an isomorphism, it follows that $b_s=0$.
Consequently,
\begin{align}
\psi^{-1}(\cA\langle r\rangle)&=\tUi\langle r\rangle,
&\cA\langle r\rangle*\cA\langle s\rangle
&\subseteq\cA\langle r+s\rangle.
\label{eq:graded-star}
\end{align}
Hence,
\begin{align*}
x*y\in\cA\langle m+n\rangle,\qquad\cA^{[0]}\langle r\rangle=\cR_{\bullet,i}^{-,r}[r]\cR_i\cH_{\bullet,i}.
\end{align*}
Since $\deg_{\bullet,i}^-\tk_\nu=2d_{\bullet,i}(\nu)$ and $\deg_{\bullet,i}^-(\cR_i\cH_{\bullet,i})=0$, comparison of homogeneous components in \eqref{eq:pJ} yields
\begin{align*}
z_{\nu}
\in\cA^{[0]}\langle m+n-2d_{\bullet,i}(\nu)\rangle
=\cR_{\bullet,i}^{-,r}[m+n-2d_{\bullet,i}(\nu)]
\cR_i\cH_{\bullet,i}.
\end{align*}
\end{proof}

\begin{lemma}\label{lem:specialization}
Let $h\in\cA^{[0]}$ and set $c=\left.\rho_t(\psi^{-1}(h))\right|_{t=0}$. Then
$c\in\tU^\poly$ and $\psi(c)=h$.
\end{lemma}
\begin{proof}
 Equation~\eqref{eq:polynomial and constant-z} gives $c\in\tU^\poly$. Together with \eqref{eq:rho-psi} and $\rho_t(\psi^{-1}(h))\in\tU^\poly[t]$, it follows that 
 \begin{equation*}
  h=\psi\big(\rho_t(\psi^{-1}(h))\big)|_{t=0}=\psi\big(\rho_t(\psi^{-1}(h))|_{t=0}\big)=\psi(c)   
 \end{equation*}
\end{proof}

Theorem~\ref{thm: relative T_i conjugate} gives
\begin{align}
\tTa{i}(x)
=\tfX_i\tT_{\br_i}^{-1}(x)\tfX_i^{-1},
\qquad x\in\tUi.
\label{eq:intertwiner of iT}
\end{align}
\cite[Remark 3.10, Proposition 4.6]{WZ23} gives
\begin{align}
\tfX_i^{\pm1}
&\in 1+\prod_{m\in\Z_{>0}}\tU_{m(\alpha_i+\bw\alpha_{\tau i})}^+,
&\tTD_j^{\pm1}(\tfX_i)&=\tfX_i\quad(j\in\bI).
\label{eq:quasiK}
\end{align}
By \cite[Lemma 2.16]{KY21}, $\tTD_{\br_i}^{-1}(z)\in\tU^-$ for $z\in\cR_{\bullet,i}^{-,r}$. Hence,
\begin{align}
\tT_{\br_i}^{-1}(z)=\tTD_{\br_i}^{-1}(z),\qquad
z\in\cR_{\bullet,i}^{-,r}.
\label{eq:rescaling-negative}
\end{align}

\begin{proposition}\label{lem:Cartan}
For every $\nu\in\Z\J$,
\begin{align}
(\psi\circ\tTa{i}\circ\psi^{-1})
\bigl(\cA^{[0]}\tk_{\nu}\bigr)\subseteq\cA^{[0]}\tk_{\nu}.\label{eq:Cartan-preserve}
\end{align}
\end{proposition}

\begin{proof}
Consider the subspace
\begin{align}
\cV_{\bullet,i}=\bigoplus_{\lambda\in\Z \bIi,\,\mu\in\Z\I}
\tU^-\tU^+K_\lambda K_\mu'.
\label{eq:VJ}
\end{align}
By \eqref{eq:support-all},
\begin{align}
\tbiU\cV_{\bullet,i}+\cV_{\bullet,i}\tbiU&\subseteq\cV_{\bullet,i},
&\psi^{-1}(\cA^{[0]})
&=\psi^{-1}(\tU^-)\cH_{\bullet,i}\subseteq\cV_{\bullet,i}.
\label{eq:VJ-local}
\end{align}
For $j\in \bIi$, Proposition~\ref{prop:braid1} gives
\begin{align*}
\tTD_j^{\pm1}(\tU^-)&\subseteq\tU^-\tbiU^+\tbiU^0,
&\tTD_j^{\pm1}(\tU^+)&\subseteq\tbiU^-\tU^+\tbiU^0,\\
(\tU^-\tbiU^+\tbiU^0)(\tbiU^-\tU^+\tbiU^0)
&\subseteq\tU^-\tU^+\tbiU^0,
&s_j(\Z \bIi)&=\Z \bIi.
\end{align*}
Denote by $\widehat{\cV_{\bullet,i}}$ the completion of $\cV_{\bullet,i}$. 
The preceding argument gives
\begin{align}
\tTD_j^{\pm1}(\cV_{\bullet,i})&=\cV_{\bullet,i},
&\tPsi_{\bvs_\dm}^{\pm1}(\cV_{\bullet,i})&=\cV_{\bullet,i},\nonumber\\
\tT_{\br_i}^{-1}(\cV_{\bullet,i})&=\cV_{\bullet,i},
&\tfX_i\cV_{\bullet,i}\tfX_i^{-1}&\subseteq\widehat{\cV_{\bullet,i}}.
\label{eq:VJ-braid}
\end{align}

Let $a\in\cA^{[0]}$. Equations~\eqref{eq:intertwiner of iT} and \eqref{eq:VJ-braid} imply
\begin{align}
\tTa{i}(\psi^{-1}(a))
\in\widehat{\cV_{\bullet,i}}\cap\tU=\cV_{\bullet,i},
\label{eq:VJ-image}
\end{align}
Accordingly, write
\begin{align*}
\psi(\tTa{i}(\psi^{-1}(a)))
&=\sum_{\nu}a_{\nu}\tk_{\nu},
\qquad a_{\nu}\in\cA^{[0]},\\
\tTa{i}(\psi^{-1}(a))
&=\sum_{\nu}
\psi^{-1}(a_{\nu})\tk_{\nu}.
\end{align*}
Equations~\eqref{eq:VJ}, \eqref{eq:VJ-local}, and \eqref{eq:VJ-image} show that the $K_\lambda$-terms of $\tTa{i}(\psi^{-1}(a))$ in \eqref{eq:VJ} are supported on $\Z\bIi$. Hence $\psi^{-1}(a_{\nu})=0$ for $\nu\ne0$, proving \eqref{eq:Cartan-preserve} for $\nu=0$.

Proposition~\ref{prop:iT on bUi} gives
\begin{align*}
 \tTa{i}(\tk_{\nu})\in\cH_{\bullet,i}\tk_{\nu}.  
\end{align*}
For $a\in\cA^{[0]}$, the $\cH_\bullet$-bilinearity of $\psi$ gives
\begin{align*}
\psi\bigl(\tTa{i}(\psi^{-1}(a\tk_{\nu}))\bigr)
&=\psi\bigl(\tTa{i}(\psi^{-1}(a))
\tTa{i}(\tk_{\nu})\bigr)\\
&=\psi\bigl(\tTa{i}(\psi^{-1}(a))\bigr)
\tTa{i}(\tk_{\nu})
\in\cA^{[0]}\tk_{\nu}.
\end{align*}
This proves \eqref{eq:Cartan-preserve}. 
\end{proof}

\subsection{The relative braid restriction}
\begin{lemma}\label{lem:conjugation}
Let $v\in\sigma(\cR_{\bullet,i}^{-,r})$. If
$c=\tfX_i v\tfX_i^{-1}\in\tU^\poly$, then $\psi(c)=v$.
\end{lemma}

\begin{proof}
By \cite[Lemma 10.1]{J96} and \eqref{eq:EyyE}, one has
\begin{align}
\partial_k^L \partial_l^R&=\partial_l^R\partial_k^L,\quad k,l\in\I,\label{eq:commute of derviation}\\
E_j v - v E_j &= \frac{K_j \partial_j^L(v)-\partial_j^R(v)K_j'}{q_j-q_j^{-1}}=\frac{-\partial_j^R(v)K_j'}{q_j-q_j^{-1}},\quad j\in\bIi, v\in \sigma(\cR_{\bullet,i}^{-,r}).
\label{eq:commutator}
\end{align}
By \eqref{eq:radicals} and \eqref{eq:commute of derviation}, $\partial_j^R(v)\in\sigma(\cR_{\bullet,i}^{-,r})$ for $j\in\bIi$.

 Fix $e\in(\tfX_i)_{\alpha}$ and $g\in(\tfX_i^{-1})_{\eta}$, where $\alpha,\eta\in\N(\alpha_i+\bw\alpha_{\tau i})$. Repeated application of \eqref{eq:commutator} to $ev$ yields
 \begin{align}
evg\in\sum_{\substack{\beta,\gamma\in\N \bIi\\
\beta+\gamma=\alpha+\eta}}
\sigma(\cR_{\bullet,i}^{-,r})(\tbiU^+)_{\gamma}K_\beta'.
\label{eq:normal-ordering}
\end{align}
If $\alpha+\eta\neq0$, at least one of $\beta$ and $\gamma$ is supported on $\{i,\tau i\}$. Moreover,
\begin{align*}
 \tfX_i v\tfX_i^{-1}-v\in\Big(\prod_{\substack{\alpha,\eta\in\N(\alpha_i+\bw\alpha_{\tau i})\\ \alpha+\eta\in\Z_{>0}(\alpha_i+\bw\alpha_{\tau i})}}\prod_{\substack{\beta,\gamma\in\N \bIi\\
\beta+\gamma=\alpha+\eta}}\sigma(\cR_{\bullet,i}^{-,r})(\tbiU^+)_{\gamma}K_\beta'\Big)\cap\tU^\poly.   
\end{align*}
Thus,
\begin{equation*}
  \tfX_i v\tfX_i^{-1}-v\in\Big(\sum_{\substack{\alpha,\eta\in\N(\alpha_i+\bw\alpha_{\tau i})\\ \alpha+\eta\in\Z_{>0}(\alpha_i+\bw\alpha_{\tau i})}}\sum_{\substack{\beta,\gamma\in\N \bIi\\
\beta+\gamma=\alpha+\eta}}\sigma(\cR_{\bullet,i}^{-,r})(\tbiU^+)_{\gamma}K_\beta'\Big)\cap\tU^\poly.  
\end{equation*}
 Hence $\psi(\tfX_i v\tfX_i^{-1}-v)=0$.
\end{proof}

\begin{theorem}\label{thm:restriction}
For $z\in\cR_{\bullet,i}^{-,r}$, one has
\begin{align}
\psi\circ\tTa{i}\circ\psi^{-1}(z)=\tTD_{\br_i}^{-1}(z)\in\tU^-.
\label{eq:restriction}
\end{align}
\end{theorem}

\begin{proof}
Let $z\in\cR_{\bullet,i}^{-,r}$ and set
\begin{align}
h=\psi\bigl(\tTa{i}(\psi^{-1}(z))\bigr).
\label{eq:h-definition}
\end{align}
By Proposition \ref{lem:Cartan}, $h\in\cA^{[0]}$.
Since $\tTa{i}(\psi^{-1}(z))\in\tUi$, \eqref{eq:intertwiner of iT} gives
\begin{align}
\psi^{-1}(h)
=\tfX_i\tT_{\br_i}^{-1}(\psi^{-1}(z))\tfX_i^{-1}.
\label{eq:lifted-intertwiner}
\end{align}
Lemma~\ref{lem:twisting Psi}, Proposition~\ref{prop:braid1}, and \eqref{eq:quasiK} imply
\begin{align}
\rho_t(\tfX_i^{\pm1})&=\tfX_i^{\pm1},
&\rho_t\tT_{j}^{-1}&=\tT_{j}^{-1}\rho_t,\quad j\in\bIi.
\label{eq:rho-local}
\end{align}

Equation~\eqref{eq:polynomial and constant-z} shows that $\rho_t(\psi^{-1}(h)),\rho_t(\psi^{-1}(z))\in\tU^\poly[t]$. Apply $\rho_t$ to \eqref{eq:lifted-intertwiner} and then set $t=0$. Equations~\eqref{eq:polynomial and constant-z} and \eqref{eq:rho-local} give
\begin{align}
\left.\rho_t(\psi^{-1}(h))\right|_{t=0}
&=\tfX_i\tT_{\br_i}^{-1}
\left(\left.\rho_t(\psi^{-1}(z))\right|_{t=0}\right)
\tfX_i^{-1}\nonumber\\
&=\tfX_i\tT_{\br_i}^{-1}(z)\tfX_i^{-1}
=\tfX_i\tTD_{\br_i}^{-1}(z)\tfX_i^{-1}\in\tU^\poly.
\label{eq:specialized-intertwiner}
\end{align}

Define
\begin{align*}
c=\tTD_{\bw}^{-1}
\left(\left.\rho_t(\psi^{-1}(h))\right|_{t=0}\right).
\end{align*}
Lemma~\ref{lem:specialization} gives
\begin{align}
c\in\tU^\poly,
\qquad\psi(c)=\tTD_{\bw}^{-1}(h).
\label{eq:c-projection}
\end{align}
Combining \eqref{eq:quasiK}, $\tTD_{\bw}^{-1}\tTD_{\br_i}^{-1}=\tTD_{\bwi}^{-1}$, and \eqref{eq:specialized-intertwiner} gives
\begin{align}
c&=\tTD_{\bw}^{-1}
\bigl(\tfX_i\tTD_{\br_i}^{-1}(z)\tfX_i^{-1}\bigr)=\tfX_i\tTD_{\bwi}^{-1}(z)\tfX_i^{-1}.
\label{eq:c-conjugation}
\end{align}
By \eqref{eq:radicals},
$v=\tTD_{\bwi}^{-1}(z)\in\sigma(\cR_{\bullet,i}^{-,r})$.
Applying Lemma~\ref{lem:conjugation} to \eqref{eq:c-conjugation} gives
\begin{align*}
\tTD_{\bw}^{-1}(h)=\psi(c)=\tTD_{\bwi}^{-1}(z),
\qquad
h=\tTD_{\bw}\tTD_{\bwi}^{-1}(z)=\tTD_{\br_i}^{-1}(z).
\end{align*}
This proves \eqref{eq:restriction}.
\end{proof}

The following corollary extends the compatibility in \cite[Proposition 3.9]{SZ25} to symmetrizable Kac–Moody type.
\begin{corollary}
For $j\notin\bIi$, we have
\begin{align*}
\psi\circ\tTa{i}\circ\psi^{-1}(F_j)=\tTD_{\br_i}^{-1}(F_j).
\end{align*}
\end{corollary}

Let $w\in W^\circ$. Given a reduced expression $\underline{w} =\br_{i_1} \br_{i_2} \ldots \br_{i_k}$ for $w$, we shall write $\tTT'_{\underline{w},-1} =\tTT'_{i_1,-1} \tTT'_{i_2,-1} \ldots \tTT'_{i_k,-1}$.

The following corollary extends \cite[Theorem 7.13]{WZ23} to symmetrizable Kac–Moody type.
\begin{corollary}
\label{thm:fact1}
Suppose that $wi\in \wI$, for $w\in W^\circ$ and $i\in \wI$. Then $\tTT'_{\underline{w},-1}(B_i)=B_{wi}$, for {\em some} reduced expression $\underline{w}$ of $w$.    
\end{corollary}
\begin{proof}
Using Lemma \ref{lem:Tw-on-X}, we have
\begin{align*}
  \tU^-\cap\tTD_{\br_i}(\tU^-)=&\tTD_{w_\bullet}^{-1}(\tTD_{w_\bullet}(\tU^-)\cap\tTD_{w_\bullet,i}(\tU^-))\\
  =&\tTD_{w_\bullet}^{-1}(\tTD_{w_\bullet}(\tbU^-)\otimes\cR_\bullet^{-,r}\cap\tTD_{w_\bullet,i}(\tbiU^-)\otimes\cR_{\bullet,i}^{-,r})\\
  =&\tTD_{w_\bullet}^{-1}\Big(\big(\tTD_{w_\bullet}(\tbU^-)\cap\tTD_{w_\bullet,i}(\tbiU^-)\big)\otimes\cR_{\bullet,i}^{-,r}\Big)\\
  =&\tbU^-\otimes\tTD_{w_\bullet}^{-1}(\cR_{\bullet,i}^{-,r}).
\end{align*}
By \eqref{eq:intertwiner of iT}, \eqref{eq:quasiK} and Lemma~\ref{lem:Letzter-props}, $\tTD_{w_\bullet}^{\pm1}\circ\psi\circ\tTa{i}\circ\psi^{-1}=\psi\circ\tTa{i}\circ\psi^{-1}\circ\tTD_{w_\bullet}^{\pm1}$.
So for $x\in\tU^-\cap\tTD_{\br_i}(\tU^-)$, using Theorem \ref{thm:restriction}, we have
\begin{equation}\label{eq:general of thm:res}
    \psi\circ\tTa{i}\circ\psi^{-1}(x)=\tTD_{\br_i}^{-1}(x)\in\tU^-.
\end{equation}

Given a reduced expression $\underline{w} =\br_{i_1} \br_{i_2} \ldots \br_{i_l}$ for $w$, we set $w_m=\br_{i_m} \ldots \br_{i_1}$ and $w_{0}=1$ for $0\leq m\leq l$. Since $wi\in \wI$, we have $w_m(\alpha_{wi})>0$ for $0\leq m\leq l$. Hence, we set
\begin{equation*}
    y_m=\tTD_{w_m}(F_{wi})\in\tU^-.
\end{equation*}
Note that $y_l=F_i$. Then $y_{m+1}=\tTD_{\br_{i_{m+1}}}(y_m)\in\tU^-\cap\tTD_{\br_{i_{m+1}}}(\tU^-)$. By \eqref{eq:general of thm:res}, we have
\begin{equation*}
 \psi\circ\tTa{i_{m+1}}\circ\psi^{-1}(y_{m+1}) =y_m.  
\end{equation*}
Furthermore, we have
\begin{equation*}
  \psi\circ\tTT'_{\underline{w},-1}\circ\psi^{-1}(F_i)=y_0=F_{wi}.  
\end{equation*}
Applying $\psi^{-1}$ and using $\psi^{-1}(F_i)=B_i$, we complete the proof.
\end{proof}

%%%%%%%%%%%%%%%%
%%%%%%%%%%%%%%%%
\section{\texorpdfstring{$\mathrm i$}{i}Braid group symmetries via star product}
\label{sec:global} 

In this section, we use the star product to relate the braid group symmetries of quantum groups to those of iquantum groups. 

\subsection{The opposite Letzter projection}
Recall the subalgebras $\tU^\poly$ and $\cA$, and the left ideal $\cJ_\bullet$ of $\tU^\poly$. They satisfy 
\begin{align*}
 \tU^\poly=\cA\oplus\cJ_\bullet,\qquad
\psi:\tU^\poly\longrightarrow\cA.   
\end{align*}
Applying $\sigma$ to the construction in Section~\ref{sec:letzter map} gives
\begin{align}
&\sigma(\tU^\poly)=\cA\oplus\sigma(\cJ_\bullet),\qquad\sigma(\tU^\poly)\cong\sigma(\cR_{\bullet}^{+,r})\otimes \K[K_i'^{-1}\,|\,i\in \I_\circ]\otimes \cA,\label{eq:opposite-Upoly}\\
&\sigma(\cJ_\bullet)=\Bigl(\sigma(\cR_{\bullet,+}^{+,r})
+\sum_{j\in\wI}\K K_j'^{-1}\Bigr)\sigma(\tU^\poly)= \cH_\bullet\Span_\K\{E_jK_j'^{-1},K_j'^{-1}\mid j\in\wI\}\sigma(\tU^\poly).\label{eq:right idea}
\end{align}
By \eqref{eq:crucial-deg}, $\sigma(\tU^\poly)$ is a graded subalgebra for $\deg_{\bullet,i}^+$ and $\sigma(\cJ_\bullet)$ is a graded right ideal of $\sigma(\tU^\poly)$.
Define
\begin{equation}
\psi^\sigma:=\sigma\circ\psi\circ\sigma:\sigma(\tU^\poly)\rightarrow\cA.\label{eq:opposite-projection}
\end{equation}
Thus $\psi^\sigma$ is the $\K$-linear projection associated with the direct-sum decomposition \eqref{eq:opposite-Upoly}.

For the automorphism $\rho_t$ of $\tU\otimes_\K\K[t,t^{-1}]$,
\begin{equation}\label{eq:sigma and specialization}
  \rho_t\circ\sigma=\sigma\circ\rho_t,\qquad  \rho_t(\tfX)\big|_{t=0}=\tfX_i. 
\end{equation}

\begin{lemma}\label{lem:opposite-braid}
For $a\in\cA^{[0]}$, one has
\begin{align}
\tTD_{\bw}^{-1}\circ\psi\circ\tTa{i}\circ\psi^{-1}(a)
=\psi^\sigma\left(
\tT_{\bwi}^{-1}\left(
\left.\rho_t(\psi^{-1}(a))\right|_{t=0}\right)\right).
\label{eq:opposite-braid}
\end{align}
\end{lemma}

\begin{proof}
\cite[Corollary 4.7]{KY26} and Lemma~\ref{lem:anti-involut-QG} give
\begin{align*}
  \sigma^\imath=\psi^{-1}\circ\sigma\circ\psi,\qquad \sigma^\imath(x)\tfX&=\tfX\sigma(x), \quad x\in\tUi. 
\end{align*}
For $h\in\cA^{[0]}$,
\begin{equation}
  \sigma\circ\psi^{-1}\circ\sigma(h)=\tfX^{-1}\psi^{-1}(h)\tfX.\label{eq:sigma-lift}  
\end{equation}

Set $b_0=\left.\rho_t(\psi^{-1}(a))\right|_{t=0}$ and
$h=\tTD_{\bw}^{-1}\big(\psi\circ\tTa{i}\circ\psi^{-1}(a)\big)$. By \eqref{eq:polynomial and constant-z} and \eqref{eq:Cartan-preserve}, $b_0$ is well defined and $h\in\cA^{[0]}$.
Lemma~\ref{lem:Letzter-props}, equations~\eqref{eq:quasiK}, \eqref{eq:polynomial and constant-z}, \eqref{eq:rho-local}, and $\tTD_{\bw}^{-1}\tT_{\br_i}^{-1}=\tT_{\bwi}^{-1}$ give
\begin{align*}
\left.\rho_t(\psi^{-1}(h))\right|_{t=0}
&=\tfX_i\tT_{\bwi}^{-1}(b_0)\tfX_i^{-1}.
\end{align*}
Equations~\eqref{eq:polynomial and constant-z} and \eqref{eq:sigma and specialization} imply $\rho_t\bigl(\sigma\circ\psi^{-1}\circ\sigma(h)\bigr)\in\sigma(\tU^\poly)[t]$.
Setting $t=0$ and using \eqref{eq:rho-local} and \eqref{eq:sigma-lift} gives
\begin{align}
\left.\rho_t\bigl(\sigma\circ\psi^{-1}\circ\sigma(h)\bigr)\right|_{t=0}
&=\tT_{\bwi}^{-1}(b_0)\in\sigma(\tU^\poly).
\label{eq:opposite-lift}
\end{align}
Applying $\psi^\sigma$ to \eqref{eq:opposite-lift}, together with \eqref{eq:sigma and specialization}, proves \eqref{eq:opposite-braid}.
\end{proof}

The triangular decompositions of $\tU$ define the projections
\begin{align}
    \pi^{\le0}:&\tU=\tU^+\tU^-\tU^0\longrightarrow \tU^-\tU^0,\label{eq:proj-nega}\\
    \pi^{\ge0}:&\tU=\tU^-\tU^+\tU^0\longrightarrow \tU^+\tU^0.\label{eq:proj-posi}
\end{align}

\begin{lemma}\label{lem:PBW}
The following properties hold.
\begin{enumerate}
\item $\pi^{\le0}(\sigma(\tU^\poly))\subseteq\sigma(\tU^\poly)$.
\item For $z\in\tU^-$ and
$v\in\sigma(\tU^\poly)\cap\tU^-\tU^0$,
\begin{align}
\psi^\sigma(zv)=z\psi^\sigma(v).\label{eq:negative-Cartan}
\end{align}
\item For homogeneous $b\in\sigma(\cR_i)\cap\ker\eps$,
\begin{align}
\tT_{\bwi}^{-1}(b)\in\sigma(\cJ_\bullet).
\label{eq:positive-kernel}
\end{align}
\end{enumerate}
\end{lemma}
\begin{proof}
    The identity $\cA=\tbU^+\tU^-\tU^{\imath0}$ and \eqref{eq:opposite-Upoly} prove (1). Moreover,
    \begin{align*}
        \tU^-\cA\subseteq\cA,\quad \tU^-\big(\sigma(\cJ_\bullet)\cap\tU^-\tU^0\big)\subseteq\sigma(\cJ_\bullet).
    \end{align*}
    This proves (2).

    Lemma~\ref{lem:Tw-on-X} gives
    \begin{align*}
        \tTD_{\bwi}^{-1}(F_j)=-E_{\tau_{\bullet,i} j} K_{\tau_{\bullet,i}j}'^{-1},\quad j\in\bIi.
    \end{align*}
    This proves (3) by \eqref{eq:generator of radicals} and \eqref{eq:right idea}.
\end{proof}

\subsection{The product formula}

\begin{lemma}\label{lem:c-definition2}
 For $a\in(\tbiU^-)_{-\alpha}$, there exists an element $c(a)\in\tbiU^+$ such that
 \begin{align}
     \pi^{\ge0}(\psi^{-1}(a))=c(a)K_{\alpha}',\quad c(1)=1.\label{eq:c-definition2}
 \end{align}
\end{lemma}
\begin{proof}
    Recall the proof of Lemma~\ref{lem:support}, where
    \begin{align*}
      &\Theta=\sum_{\mu,r_\mu}u^-_{\mu,r_\mu}\otimes u^+_{\mu,r_\mu},
\qquad\langle u^-_{\mu,r_\mu},u^+_{\nu,r_\nu}\rangle=\delta_{\mu,\nu}\delta_{r_\mu ,r_\nu},\\
&(\psi^{-1}\otimes\id)(\Theta)=\Theta C,\qquad
C=\bDelta(\tfX)(\tfX^{-1}\otimes1)=\sum_{\beta\in\N\I}C_\beta,\\
&C_0=1\otimes 1,\qquad C_\beta=\sum_{s_\beta} c_{\beta,s_\beta}K_\beta'\otimes v^+_{\beta,s_\beta}\in\widehat{\tU^+}K_\beta'
\widehat\otimes\tU^+_\beta.
    \end{align*}
    Equation~\eqref{eq:inverse-expansion} gives
    \begin{equation*}
        \psi^{-1}(a)=a+\sum_{\substack{\beta\ne0,\,\mu\\r_\mu,s_\beta}}
u^-_{\mu,r_\mu}c_{\beta,s_\beta}
\langle a,u^+_{\mu,r_\mu}v^+_{\beta,s_\beta}\rangle K_\beta'.
    \end{equation*}
    Under $\pi^{\ge0}$, a term of $\psi^{-1}(a)$ vanishes unless $\mu=0$ and $\beta=\alpha$. Hence
    \begin{equation*}
     \pi^{\ge0}(\psi^{-1}(a))=\sum_{s_{\alpha}}c_{\alpha,s_\alpha}
\langle a,v^+_{\alpha,s_\alpha}\rangle K_\alpha'=c(a)K_{\alpha}'.   
    \end{equation*}
By Remark~\ref{rem:local Letzter and Letzter} and Lemma~\ref{lem:support}, 
\begin{equation*}
    \psi^{-1}(a)\in\sum_{\beta\in\N\bIi}\tbiU^-\tbiU^+K_\beta'.
\end{equation*}
Hence $c(a)\in\tbiU^+$.
\end{proof}

Lemma~\ref{lem:c-definition2} also gives
\begin{equation}\label{eq:c-def of a}
    c(a)=\sum_{s_{\alpha}}c_{\alpha,s_\alpha}
\langle a,v^+_{\alpha,s_\alpha}\rangle\in\tbiU^+,\qquad a\in(\tbiU^-)_{-\alpha}.
\end{equation}

\begin{proposition}\label{lem:mixed-lift}
Let $u\in(\tbiU^-)_{-\lambda}$ and
$x\in(\cR_{\bullet,i}^{-,r})_{-\eta}$. Write
\begin{align}
\Delta(u)=\sum_m a_m\otimes b_mK_{\alpha_m}',\qquad a_m\in(\tbiU^-)_{-\alpha_m},
b_m\in(\tbiU^-)_{-(\lambda-\alpha_m)}.\label{eq:coproduct}
\end{align}
Then
\begin{align}
\psi^{-1}(u)
&=\sum_m b_m c(a_m)K_{\alpha_m}',\label{eq:local-lift}\\
\left.\rho_t(\psi^{-1}(ux))\right|_{t=0}
&=\sum_m q^{(\alpha_m,\eta)}
b_mx c(a_m)K_{\alpha_m}'.\label{eq:mixed-lift}
\end{align}
\end{proposition}

\begin{proof}
Equation~\eqref{eq:inverse-expansion} gives
\begin{align}
\psi^{-1}(u)=\sum_{\substack{\beta,\,\mu\\r_\mu,s_\beta}}
u^-_{\mu,r_\mu}c_{\beta,s_\beta}
\langle u,u^+_{\mu,r_\mu}v^+_{\beta,s_\beta}\rangle K_\beta',\label{eq:inverse-1}\\
\psi^{-1}(ux)=\sum_{\substack{\beta,\,\mu\\r_\mu,s_\beta}}
u^-_{\mu,r_\mu}c_{\beta,s_\beta}
\langle ux,u^+_{\mu,r_\mu}v^+_{\beta,s_\beta}\rangle K_\beta'.\label{eq:inverse-2}
\end{align}
Equation~\eqref{eq:form-kow} gives 
\begin{align}
 \langle u,u^+_{\mu,r_\mu}v^+_{\beta,s_\beta}\rangle=& \sum_m \langle a_m,v^+_{\beta,s_\beta}\rangle\langle b_mK_{\alpha_m}',u^+_{\mu,r_\mu}\rangle= \sum_m \langle a_m,v^+_{\beta,s_\beta}\rangle\langle b_m,u^+_{\mu,r_\mu}\rangle,\label{eq:inverse of pair-1}\\
 \langle ux,u^+_{\mu,r_\mu}v^+_{\beta,s_\beta}\rangle=&\langle u_{(1)}x_{(1)},v^+_{\beta,s_\beta}\rangle\langle u_{(2)}x_{(2)},u^+_{\mu,r_\mu}\rangle\label{eq:inverse of pair-2}.
\end{align}

In \eqref{eq:inverse of pair-1}, a term vanishes unless $\beta=\alpha_m$ for some $m$. Equation~\eqref{eq:c-def of a} therefore gives
\begin{align*}
 \psi^{-1}(u)=&\sum_{m}\sum_{s_{\alpha_m}}\sum_{\mu,r_\mu}
u^-_{\mu,r_\mu}c_{\alpha_m,s_{\alpha_m}}\langle a_m,v^+_{\alpha_m,s_{\alpha_m}}\rangle\langle b_m,u^+_{\mu,r_\mu}\rangle K_{\alpha_m}' \\
=&\sum_{m}\sum_{s_{\alpha_m}}b_mc_{\alpha_m,s_{\alpha_m}}\langle a_m,v^+_{\alpha_m,s_{\alpha_m}}\rangle K_{\alpha_m}'\\
=&\sum_m b_m c(a_m)K_{\alpha_m}'.
\end{align*}

Substituting \eqref{eq:inverse of pair-2} into \eqref{eq:inverse-2} gives
\begin{align}
 \psi^{-1}(ux)=\sum_{\substack{\beta,\,\mu\\r_\mu,s_\beta}}
u^-_{\mu,r_\mu}c_{\beta,s_\beta}
\langle u_{(1)}x_{(1)},v^+_{\beta,s_\beta}\rangle\langle u_{(2)}x_{(2)},u^+_{\mu,r_\mu}\rangle K_\beta'.\label{eq:inverse-3}   
\end{align}
Since $\cR_{\bullet,i}^{-,r}$ is a right coideal subalgebra, the terms in \eqref{eq:inverse-3} with $x_{(1)}\ne1$ have $\beta\in\N\I-\N\bIi$. By \eqref{eq:polynomial and constant-z}, $\rho_t(\psi^{-1}(ux))\in\tU^\poly[t]$. Hence
\begin{align*}
 \left.\rho_t(\psi^{-1}(ux))\right|_{t=0}=& \rho_t\Big(\sum_{m}\sum_{s_{\alpha_m}}\sum_{\mu,r_\mu}
u^-_{\mu,r_\mu}c_{\alpha_m,s_{\alpha_m}}\langle a_m,v^+_{\alpha_m,s_{\alpha_m}}\rangle\langle b_mK_{\alpha_m}'x,u^+_{\mu,r_\mu}\rangle K_{\alpha_m}'\Big) \Big|_{t=0} \\
=&\sum_mq^{(\alpha_m,\eta)} b_mx \left.\rho_t\big(c(a_m)\big)\right|_{t=0}K_{\alpha_m}'. 
\end{align*}
By \eqref{eq:c-def of a}, $c(a_m)\in\tbiU^+$. Hence
\begin{align*}
    \left.\rho_t(\psi^{-1}(ux))\right|_{t=0}=\sum_m q^{(\alpha_m,\eta)}b_mx c(a_m)K_{\alpha_m}'.
\end{align*}
\end{proof}

\begin{lemma}\label{lem:domain}
For $a\in(\tbiU^-)_{-\alpha}$, one has
\begin{align}
\tT_{\bwi}^{-1}\bigl(c(a)K_\alpha'\bigr)
\in\sigma(\tU^\poly)\cap\tbiU^-\tbiU^0.
\label{eq:coefficient-domain}
\end{align}
\end{lemma}

\begin{proof}
Since $\psi^{-1}(a)\in\tbiU$, \eqref{eq:opposite-lift} gives
$\tT_{\bwi}^{-1}(\psi^{-1}(a))\in\sigma(\tU^\poly)$. Moreover, 
\begin{equation*}
  \psi^{-1}(a)\in\sum_{\beta\in\N\bIi}\tbiU^-\tbiU^+K_\beta',\qquad \pi^{\ge0}(\psi^{-1}(a))=c(a)K_{\alpha}' .
\end{equation*}
Applying $\bwi$ to Lemma~\ref{lem:Tw-on-X} and using Lemma~\ref{lem:c-definition2} yields
\begin{equation*}
    \pi^{\le0}(\tT_{\bwi}^{-1}(\psi^{-1}(a)))=\tT_{\bwi}^{-1}(\pi^{\ge0}(\psi^{-1}(a)))=\tT_{\bwi}^{-1}\bigl(c(a)K_\alpha'\bigr).
\end{equation*}
Thus, Lemma \ref{lem:PBW}(1) gives $\tT_{\bwi}^{-1}\bigl(c(a)K_\alpha'\bigr)
\in\sigma(\tU^\poly)\cap\tbiU^-\tbiU^0$.
\end{proof}

\begin{theorem}\label{thm:mixed}
For $u\in\bigl(\tTD_{\bw}^{-1}(\cR_i)\bigr)_{-\mu}$ and $x\in(\cR_{\bullet,i}^{-,r})_{-\nu}$, the following identity holds:
\begin{align}
(\psi\circ\tTa{i}\circ\psi^{-1})(ux)
&=q^{(\mu,\nu)}\tTD_{\br_i}^{-1}(x)
(\psi\circ\tTa{i}\circ\psi^{-1})(u).
\label{eq:main}
\end{align}
\end{theorem}

\begin{proof}
We use the following elementary identity: if $v\in\tU^-_{-\mu}$ and $\Delta(v)=\sum_mx_m\otimes y_mK_{\beta_m}'$ for $x_m\in\tU^-_{-\beta_{m}}$ and $y_m\in\tU^-_{-(\mu-\beta_{m})}$, then $\Delta(\sigma(v))=\sum_m\sigma(y_m)\otimes\sigma(x_m)K_{\mu-\beta_m}'$.

Since $\tTD_{\bw}^{-1}(\cR_i)=\sigma(\cR_i)$ and $\cR_i$ is a right coideal subalgebra, write
\begin{equation*}
    \Delta(u)=\sum_m a_m\otimes b_mK_{\alpha_m}',\qquad a_m\in(\tbiU^-)_{-\alpha_m},
b_m\in(\sigma(\cR_i))_{-(\mu-\alpha_m)}.
\end{equation*}
The component with $b_m=1$ is $u\otimes K_\mu'$.
Let
\begin{align*}
z=\tTD_{\bwi}^{-1}(x)\in\sigma(\cR_{\bullet,i}^{-,r})\subset\tU^-.
\end{align*}
Since $\tPsi_{\bvs_\dm}$ fixes $\tU^-$,
$\tT_{\bwi}^{-1}(x)=z$.
Lemma~\ref{lem:opposite-braid} and Proposition \ref{lem:mixed-lift} imply
\begin{align}
\tTD_{\bw}^{-1}\circ\psi\circ\tTa{i}\circ\psi^{-1}(ux)=\sum_m q^{(\alpha_m,\nu)}
\psi^\sigma\left(
\tT_{\bwi}^{-1}(b_m)\,z\,
\tT_{\bwi}^{-1}\bigl(c(a_m)K_{\alpha_m}'\bigr)\right).
\label{eq:projected-sum}
\end{align}
If $b_m\in\ker\eps$, Lemma~\ref{lem:PBW}(3) and Lemma \ref{lem:domain} give
\begin{align*}
\tT_{\bwi}^{-1}(b_m)&\in\sigma(\cJ_\bullet),
&z\,\tT_{\bwi}^{-1}\bigl(c(a_m)K_{\alpha_m}'\bigr)
&\in\sigma(\tU^\poly).
\end{align*}
Since $\sigma(\cJ_\bullet)$ is
a right ideal of $\sigma(\tU^\poly)$, Lemma~\ref{lem:PBW}(2) gives
\begin{align}
\tTD_{\bw}^{-1}\circ\psi\circ\tTa{i}\circ\psi^{-1}(ux)
&=q^{(\mu,\nu)}z\,
\psi^\sigma\left(\tT_{\bwi}^{-1}(c(u)K_\mu')\right).
\label{eq:surviving-term}
\end{align}
Setting $x=1$ in \eqref{eq:surviving-term} gives
\begin{equation*}
  \tTD_{\bw}^{-1}\circ\psi\circ\tTa{i}\circ\psi^{-1}(u)
=\,
\psi^\sigma\left(\tT_{\bwi}^{-1}(c(u)K_\mu')\right).  
\end{equation*}
Applying $\tTD_{\bw}$ and using $\tTD_{\bw}\tTD_{\bwi}^{-1}=\tTD_{\br_i}^{-1}$ yields
\begin{equation*}
    (\psi\circ\tTa{i}\circ\psi^{-1})(ux)=q^{(\mu,\nu)}\tTD_{\br_i}^{-1}(x)
(\psi\circ\tTa{i}\circ\psi^{-1})(u).
\end{equation*}
\end{proof}

%%%%%%%%%%%%%%%%%%
%%%%%%%%%%%%%%%%%%

\subsection{An explicit formula on $\tbiU^-$}
Recall that $\mathcal{K}_\lambda$ for $\lambda\in\Z\I$ is defined in \eqref{eq:vkHla}.
\begin{proposition}\label{prop:local formula for U-}
For $u\in(\tbiU^-)_{-\mu}$ and $\mu\in\N\bIi$, one has
    \begin{align}
        \psi\circ\tTa{i}\circ\psi^{-1}(u)=\vs_{\mu,\dm}q^{(\mu,\mu+\theta\mu)/2}\tT_{\bw}^2(\tau_{\bullet,i}\circ\tau(u))\cK_{\tau_{\bullet,i}\tau \mu}^{-1}.
    \end{align}
\end{proposition}
\begin{proof}
    Recall the subalgebra $\tU_{\text{int}}$ defined in \eqref{eq:sum of bigraded subspace and int subalge}, and the endomorphisms $\wObb$, $\Xi$ defined in \eqref{eq:auto-Omega} and \eqref{eq:formula of Xi}. Since $\wObb$ is an automorphism, we define
    \begin{equation}\label{eq:def-Gamma}
        \Gamma_\bullet:=\Xi\circ\wObb^{-1}:\tU_{\text{int}}\longrightarrow\tU_{\text{int}}.
    \end{equation}
    For $\lambda\in\Z\I$, by \eqref{eq:tau-S-sig-omega retain grade} and \eqref{eq:T retain grade}, we know that $\wObb(\tU[\lambda])=\tU[\theta\lambda]$ and $\vs_{\theta\lambda,\dm}=\vs_{\lambda,\dm}^{-1}$. It follows from \eqref{eq:formula of Xi} that
    \begin{equation}\label{eq:formula of Gamma}
        \Gamma_\bullet(x)=\vs_{\lambda,\dm}q^{\frac12(\lambda,\lambda+\theta\lambda)}x\cK_{-\lambda},\qquad x\in\tU[\lambda].
    \end{equation}

    For $j\in\wI$ and $\lambda\in\Z\I$,
    \begin{equation}
        \begin{aligned}\label{eq:compute-1}
        &\Gb(F_j)=\vs_{j,\dm}q^{\frac12(\alpha_j,\alpha_j+\theta \alpha_j)}F_j\cK_j^{-1},&&\Gb(\tTD_{\bw}(E_{\tau j})K_j')=\vs_{j,\dm}q^{\frac12(\alpha_j,\alpha_j+\theta \alpha_j)}\tTD_{\bw}(E_{\tau j})K_j'\cK_j^{-1},\\
        &\Gb(B_j)=- q^{(\alpha_j,\theta \alpha_j)}B_j\cK_j^{-1},&&\Gb(\cK_\lambda)=\vs^2_{\lambda,\dm}\cK_{\theta\lambda}.
    \end{aligned}
    \end{equation}   
For $j\in\bI$, 
\begin{equation}\label{eq:compute-2}
    \Gb(F_j)=q_j^2F_jK_j^{-1}K_j',\quad \Gb(E_j)=q_j^2E_jK_jK_j'^{-1},\quad
    \Gb(K_j)=K_j,\quad \Gb(K_j')=K_j'.
\end{equation}
By \eqref{eq:sum of bigraded subspace and int subalge} and \eqref{eq:Ui and sigmaUi in Uint}, we know $\cA, \tUi\subseteq\tU_{\text{int}}$.
Thus, $\Gb(\cA)\subseteq\cA$ and $\Gb(\tUi)\subseteq\tUi$. Using Lemma~\ref{lem:Tw-on-X}, we obtain $\Gb|_{\tbU}=\tTD_{\bw}^{-2}|_{\tbU}$.

Recall that the algebra $\tU$ is bigraded by \eqref{eq:ZI-grade-1} and \eqref{eq:ZI-grade-2}. Since $\psi:\tU^\poly\rightarrow\cA$ is the projection onto a bigraded summand and $\tU^\poly,\cA,\cJ_\bullet$ are bigraded, $\psi$ is a bigraded map. Thus, for $\lambda\in\Z\I$ and $x\in\tU^\poly\cap\tU[\lambda]$, $\psi(x)\in\cA\cap\tU[\lambda]$. Since $\cK_{-\lambda}\in\tU^{\imath 0}$, by Lemma~\ref{lem:Letzter-props}, we have
\begin{equation*}
    \psi(\Gb(x))=\vs_{\lambda,\dm}q^{\frac12(\lambda,\lambda+\theta\lambda)}\psi(x\cK_{-\lambda})=\vs_{\lambda,\dm}q^{\frac12(\lambda,\lambda+\theta\lambda)}
\psi(x)\cK_{-\lambda}=\Gb(\psi(x)).
\end{equation*}
 Hence,
 \begin{equation}\label{eq:compatible of cG and psi}
     \psi\circ\Gb(x)=\Gb\circ\psi(x),\qquad x\in\tU^\poly\cap\tU_{\text{int}}.
 \end{equation}
 
Let $\tbiU^\imath$ be the iquantum group  associated with $(\bIi=\{i,\tau i\}\cup\bI,\tau|_{\bIi})$, viewed as a subalgebra of $\tUi$. By \eqref{eq:compute-1}, \eqref{eq:compute-2} and Proposition \ref{prop:iT on bUi}, a direct check on generators gives
 \begin{equation}\label{eq:local-T}
     \tTT_{i,-1}'\big|_{\tbiU^\imath}=\tT_{\bw}^2\circ\tau_{\bullet,i}\circ\tau\circ\Gb\big|_{\tbiU^\imath}.
 \end{equation}
 For $u\in(\tbiU^-)_{-\mu}$, \eqref{eq:compatible of cG and psi} gives
 \begin{align*}
     \psi\circ\tTa{i}\circ\psi^{-1}(u)=&\psi\circ\tT_{\bw}^2\circ\tau_{\bullet,i}\circ\tau\circ\Gb\circ\psi^{-1}(u)\\
     =&\tT_{\bw}^2\circ\tau_{\bullet,i}\circ\tau\circ\psi\circ\Gb\circ\psi^{-1}(u)\\
     =&\tT_{\bw}^2\circ\tau_{\bullet,i}\circ\tau\circ\Gb(u)\\
     =&\vs_{\mu,\dm}q^{\frac12(\mu,\mu+\theta\mu)}\tT_{\bw}^2(\tau_{\bullet,i}\circ\tau(u))\cK_{\tau_{\bullet,i}\tau \mu}^{-1}.
 \end{align*}
\end{proof}

\iffalse
\begin{remark}
  Equations~\eqref{eq:compatible of cG and psi} and \eqref{eq:local-T} give a formula for $\psi\circ\tTa{i}\circ\psi^{-1}$ on $\tbiU\cap\cA$.
  More precisely, for a bihomogeneous element $u\in\tbiU\cap\cA$ with 
  \begin{equation*}
  \wt(u)=\gamma,\qquad \deg_\I(u)=\xi=\beta+(1-\theta)\lambda,\quad \beta\in\Z\bI,\lambda\in\Z\bIi,   
 \end{equation*}
 we have
 \begin{align*}
     \psi\circ\tTa{i}\circ\psi^{-1}(u)=&\tT_{\bw}^2\circ\tau_{\bullet,i}\circ\tau\circ\cG(u)\\
     =&\vs_{\xi-\gamma,\dm}q^{(\gamma,\gamma+\theta \gamma)/2}\tTD_{\bw}^2(\tau_{\bullet,i}\circ\tau(u))\cK_{\tau_{\bullet,i}\tau\big(\gamma+(\theta-1)\lambda\big)}.
 \end{align*}
\end{remark}
\fi
The decomposition $\cA\cong\cH_\bullet\otimes\sigma(\cR_i)\otimes\cR_{\bullet,i}^{-,r}$, together with Propositions~\ref{prop:iT on bUi} and \ref{prop:local formula for U-} and Theorems~\ref{thm:restriction} and \ref{thm:mixed}, yields the following result.
\begin{theorem}\label{thm:full-action}
  For $h\in\cH_\bullet$, $u\in\bigl(\tTD_{\bw}^{-1}(\cR_i)\bigr)_{-\mu}$, and $x\in(\cR_{\bullet,i}^{-,r})_{-\nu}$, one has
\begin{align}
\begin{split}
(\psi\circ\tTa{i}\circ\psi^{-1})(hux)=&\vs_{\mu,\dm}q^{\frac12(\mu,\mu+\theta\mu+2\nu)}\tT_{\br_i}^{-1}(h)\cdot\tTD_{\br_i}^{-1}(x)\cdot\tT_{\bw}^2(\tau_{\bullet,i}(\tau(u)))\cK_{\tau_{\bullet,i}\tau \mu}^{-1}.
\end{split}
\end{align}  
\end{theorem}

\iffalse
\begin{theorem}\label{thm:full-action}
  For $\nu\in\Z\wI$, $y\in\tbU$, $u\in\bigl(\tTD_{\bw}^{-1}(\cR_i)\bigr)_{-\mu}$, and $x\in(\cR_{\bullet,i}^{-,r})_{-\nu}$, one has
\begin{align}
\begin{split}
(\psi\circ\tTa{i}\circ\psi^{-1})(yux\tk_{\nu})=&\vs_{\mu+\nu-\br_i\nu,\dm}q^{(\mu,\mu+\theta\mu+2\nu)}\tau_{\bullet,i}(\tau(y))\cdot\tTD_{\br_i}^{-1}(x)\\
&\times\tTD_{\bw}^2(\tau_{\bullet,i}(\tau(u)))\cdot K_{\br_i\nu-\tau_{\bullet,i}\tau \mu}
K'_{\br_i\tau\nu+\theta(\tau_{\bullet,i}\tau \mu)}.
\end{split}
\end{align}  
\end{theorem}
\fi

\end{document}